\documentclass[12pt]{article}\topmargin=0cm
\usepackage{graphicx}   %
\usepackage{enumerate}  %
\usepackage{latexsym}   %
\usepackage{amsmath}    %
\usepackage{amsbsy}     %
\usepackage{amsfonts}   %
\usepackage{mathrsfs}       
\usepackage{amssymb}    
\usepackage{array}      %
\usepackage{theorem}
\usepackage{color}

\definecolor{winered}{rgb}{0.8,0,0}

\usepackage[colorlinks=true]{hyperref}
\hypersetup{urlcolor=blue, citecolor=red}

 \newtheorem{thm}{Theorem}[section]
 \newtheorem{mainTh}{Main Theorem}
 
 \newtheorem{lem}[thm]{Lemma}
 
 \theoremstyle{definition}
 \newtheorem{defn}{Definition}[section]
 \newtheorem{rem}{Remark}[section]
 
 \numberwithin{equation}{section}

\def\R{{\mathbb R}}  
\def\N{{\mathbb{N}}}

\def\L{{\mathcal L}}
\def\M{{\mathcal M}}

\def\ds{\displaystyle}

\def\laplace{\varDelta}

\def\laplace{{\varDelta}}

\begin{document}
\label{page:t}

\begin{center}
  {\large \bf {\bf
	Energy dissipative solutions to \\[0.75ex] the Kobayashi-Warren-Carter system\footnotemark[1]
  }}
\end{center}
\vspace{7ex}
\begin{center}
{\sc Salvador Moll}\footnotemark[2]\\
Department d'An$\grave{{\rm a}}$lisi Matem$\grave{{\rm a}}$tica, Universitat de Val$\grave{{\rm e}}$ncia\\
C/Dr. Moliner, 50, Burjassot, Spain\\
{\ttfamily j.salvador.moll@uv.es}
\vspace{2ex}

{\sc Ken Shirakawa}\footnotemark[3]\\
Department of Mathematics, Faculty of Education, Chiba University\\
1-33, Yayoi-cho, Inage-ku, Chiba, 263-8522, Japan\\
{\ttfamily sirakawa@faculty.chiba-u.jp}
\vspace{2ex}

              {\sc Hiroshi Watanabe}\footnotemark[4]\\
Department of General Education, Salesian Polytechnic\\
4-6-8, Oyamagaoka, Machida-city, Tokyo, 194-0215, Japan\\
{\ttfamily h-watanabe@salesio-sp.ac.jp}
\end{center}
\vspace{10ex}

\noindent
{\bf Abstract.}
In this paper we study a variational system of two parabolic PDEs, called the Kobayashi-Warren-Carter system, which models the grain boundary motion in a polycrystal. The focus of the study is the existence of solutions to this system which dissipate the associated energy functional. We obtain existence of this type of solutions via a suitable approximation of the energy functional with Laplacians and an extra regularization of the weighted total variation term of the energy. As a byproduct of this result, we also prove some $\Gamma$-convergence results concerning  weighted total variations and the corresponding time-dependent cases. Finally, the regularity obtained for the solutions together with the energy dissipation property, permits us to completely characterize the $\omega$-limit set of the solutions.

\footnotetext[0]{
$\empty^*$\,AMS Subject Classification 35K87, 
	35R06, 
	35K67.  
$\empty^\dag$\,This author is supported by the Spanish MEC project MTM2012-31103.
$\empty^\ddag$\,This author is supported by Grant-in-Aid No. 24740099, JSPS.
$\empty^\S$\,This author is supported by Grant-in-Aid No. 25800086 and No. {26400138}, JSPS. 
}
\newpage
\section*{\hypertarget{Intro}{Introduction}}
\ \ \vspace{-4ex}

Let $ 1 < N \in \N $ be a fixed number,  let $ \Omega \subset \R^N $ be a bounded domain with a smooth boundary $ \partial \Omega $, and let $ \nu_{\partial \Omega} $ be the unit outer normal on $ \partial \Omega $. Let $ Q := (0, \infty) \times \Omega $ be the product space of the time-interval $ (0, \infty) $ and the spatial domain $ \Omega $. We also set $ \Sigma := (0, \infty) \times \partial \Omega $.

In this paper, we consider a system of parabolic PDEs, called \emph{Kobayashi-Warren-Carter system}. This system (denoted by (\hyperlink{S}{S})) consists of gradient flows derived from the following energy functional, called \emph{free energy}:
\begin{equation}\label{free}
[\eta, w] \in H^1(\Omega) \times BV(\Omega) \mapsto \mathscr{F}(\eta, \theta) := \frac{1}{2} \int_\Omega |\nabla \eta|^2 \, dx +\int_\Omega \hat{g}(\eta) \, dx +\int_\Omega \alpha(\eta) |D\theta|.
\end{equation}
Our system (\hyperlink{S}{S}) is formally described as follows: \medskip

\noindent
(\hypertarget{(S)}{S}):
\begin{equation}\label{1steq}
\left\{ ~ \parbox{8cm}{
$ \ds \eta_t -{\mit \Delta} \eta +g(\eta) +\alpha'(\eta) |D \theta| = 0 $~ in $ Q $,
\\[1ex]
$ \ds \nabla \eta \cdot \nu_{\partial \Omega} = 0 $~ on $ \Sigma $,
\\[1ex]
$ \eta(0, x) = \eta_0(x) $,~ $ x \in \Omega $;
} \right.
\end{equation}
\begin{equation}\label{2ndeq}
\left\{ ~ \parbox{8cm}{
$ \ds \alpha_0(\eta) \theta_t -{\rm div} \left( \alpha(\eta) \frac{D \theta}{|D \theta|} \right) = 0 $~ in $ Q $,
\\[1ex]
$ \ds \alpha(\eta) \frac{D\theta}{|D\theta|} \cdot \nu_{\partial \Omega} = 0 $~ on $ \Sigma $,
\\[1ex]
$ \theta(0, x) = \theta_0(x) $,~ $ x \in \Omega $.
} \right.
\end{equation}

The derivation of (\hyperlink{S}{S}) is based on the modelling method of a mathematical model of grain boundary motion, proposed by Kobayashi-Warren-Carter \cite{KWC1,KWC2}. In the original studies \cite{KWC1,KWC2} the spatial domain $ \Omega $ is settled as a two-dimensional domain, and the time-spatial variations of grain boundaries are reproduced by a vector field
$$
(t, x) \in Q \mapsto \eta(t, x) \left[ \rule{0pt}{10pt} \cos \theta(t, x), \sin \theta(t, x) \right],
$$
consisting of two order parameters $ \eta = \eta(t, x) $ and $ \theta = \theta(t, x) $. In the model, the dynamics of $ \eta $ and $ \theta $ are governed by the gradient flows of the \emph{free-energy.} In the context, $ \eta = \eta(t, x) $ and $ \theta = \theta(t, x) $ indicate, respectively, the orientation order and the orientation angle of the grain. In particular, $ \eta $ is supposed to satisfy the range constraint $ 0 \leq \eta \leq 1 $ on $ Q $, and the threshold values $ 1 $ and $ 0 $ indicate the completely oriented phase and the disordered phase of orientation, respectively.

 Here, the initial-boundary value problems (\ref{1steq}) and (\ref{2ndeq}) are the gradient flows of $ \mathscr{F} $ with respect to the order parameters $ \eta $ and $ \theta $, respectively. $ g = g(\eta) $ is a perturbation to realize the range constraint for $ \eta $, and $ \hat{g} $ is a nonnegative primitive of $ g $. $ \alpha_0 = \alpha_0(\eta) $ and $ \alpha = \alpha(\eta) $ are given positive-valued functions which control the mobility of grain boundaries, and $ \alpha' $ denotes the differential of $ \alpha  $. $ \eta_0 = \eta_0(x) $ and $ \theta_0 = \theta_0(x) $ are given initial data.
The integral $ \int_\Omega \alpha(\eta) |D\theta| $ in (\ref{free}) denotes the total variation of $ \theta $ with the unknown-dependent weight $ \alpha(\eta) $.

From a physical point of view, the role of the total variation $ \int_\Omega \alpha(\eta) |D\theta| $ is built in to reproduce the facet-like situations as in  crystalline structures. However, this total variation term brings down two nonstandard terms in the system (\hyperlink{(S)}{S}): $ \alpha'(\eta)|D \theta| $ and $ -{\rm div}(\alpha(\eta) \frac{D\theta}{|D \theta|}) $ which make the mathematical analysis very tough.
Concerning the system (\hyperlink{(S)}{S}), there is just one recent mathematical result, \cite{MS}, about existence of solutions. Previous studies of the Kobayashi-Warren-Carter system dealt with some simplified versions of (\hyperlink{S}{S}), such as
\begin{enumerate}
\item[--]restricted versions to one-dimensional cases of $ \Omega $ (cf. \cite{GG,GGK,KG,SW,SWY,WS13,WS14});
\item[--]relaxed versions by Laplacians (cf. \cite{IKY08,IKY09,IKY12,KY,KWC1,KWC2}).
\end{enumerate}

Here we focus on a notion of solution, named as \emph{energy-dissipative solution}, proposed in \cite{SW} for the one-dimensional case, which permits to obtain a soft smoothing effect and energy-dissipation of the energy functional (as is usual in parabolic type systems). The goal of this paper is to obtain the following results:
\begin{description}
\item[Main Theorem \ref{mainTh1}:]the existence of energy-dissipative solutions to (\hyperlink{(S)}{S}).
\item[Main Theorem \ref{mainTh2}:]the large-time behavior of energy-dissipative solutions.
\end{description}
\medskip

In order to prove these results, we use the approximations proposed in \cite{MS}; i.e. different time interpolations for the solutions of the associated elliptic systems together with a regularization of the nonstandard terms with {Laplacians}. Even if this method leaded to existence of solutions in \cite{MS}, the regularity of the approximating solutions was not enough to derive a crucial energy inequality (see Lemma \ref{lem1}) which, together with an improved regularity, permits to obtain the result of energy dissipation. In order to obtain this energy inequality and an improvement on the regularity of the approximating solutions, we have to approximate also the term $ \int_\Omega \alpha(\eta) |D\theta| $ with a ``suitable regularization'' of the {Euclidean} distance, {denoted by $|\cdot|_\nu$}. As a byproduct of our results, we prove new $\Gamma-$convergence results concerning energy functionals related to generalized weighted total variations and time-dependent weighted total variation (see Theorems \ref{lem-gamma0} and \ref{lem-gamma}).
\medskip

Here is the content of this paper. In the next Section 1, some specific notations are prepared as preliminaries. In Section 2, the Main Theorems in this paper are presented. In Section 3, we confirm the existence, uniqueness and a priori estimates for the approximating problems, which are prescribed as the time-discretization systems of some sorts of  relaxed versions of (\hyperlink{S}{S}) with Laplacians. In Section 4, we prove the $\Gamma-$convergence results as well as some related auxiliary results needed in the proofs of Main Theorems. Sections 5 and 6 are devoted to the proofs of Main Theorems 1 and 2, respectively.
\pagebreak

\hypertarget{Sec.1}{\empty}
\section{Preliminaries}
\ \ \vspace{-4ex}

We begin with some notations used throughout this paper.
\medskip

\noindent
\underline{\textit{Abstract notations. (cf. \cite[Chapter II]{Brezis})}}
For an abstract Banach space $X$, we denote by $|\cdot|_{X}$ the norm of $X$. In particular, when $X$ is a Hilbert space, we denote by $(\,\cdot\,,\cdot\,)_{X}$ the inner product of $X$. In particular, when $X = \R^{N}$, we simply denote by
\begin{equation*}
\textstyle
|\xi| := \sqrt{\xi_{1}^{2} + \cdots + \xi_{N}^{2}}  \mbox{ \ and \ } \xi \cdot \eta := \xi_{1}\eta_{1} + \cdots + \xi_{N}\eta_{N},
\end{equation*}
the Euclidean norm of $\xi \in \R^{N}$, and the usual scalar product of $\xi = ({\xi}_{1}, \ldots, \xi_{N})$ and $\eta = (\eta_{1},\ldots\eta_{N}) \in \R^{N}$, respectively.

For any proper functional $ \Psi : X \longrightarrow (0, \infty] $ on a Banach space $ X $, we denote by $ D(\Psi) $ the effective domain of $ \Psi $.

For any proper lower semi-continuous (l.s.c., in short) and convex function $\Phi$ defined on a Hilbert space $H$, we denote by $\partial \Phi$ the subdifferential of $\Phi$. The subdifferential $\partial \Phi$ corresponds to a weak differential of $\Phi$, and it is known as a maximal monotone graph in the product space $ H^2 $ $ (= H \times H) $. More precisely, for each $v_{0} \in H$, the value $\partial \Phi(v_{0})$ of the subdifferential at $v_{0}$ is defined as a set of all elements $v_{0}^{\ast} \in H$ which satisfy the following variational inequality:
\begin{equation*}
(v_{0}^{\ast}, v-v_{0})_{H} \le \Phi(v) - \Phi(v_{0})\ \ \mbox{for any}\ v \in D(\Phi).
\end{equation*}
The set $D(\partial \Phi) := \{ z \in H \ |\ \partial\Phi(z) \neq \emptyset\}$ is called the domain of $\partial\Phi$. We often use the notation ``$ [v_{0},v_{0}^{\ast}] \in \partial\Phi_{0}$ in $ H^2 $\,'', to mean that $`` v_{0}^{\ast} \in \partial\Phi(v_{0})$ in $H$ for $v_{0} \in D(\partial\Phi)"$, by identifying the operator $\partial\Phi$ with its graph in $H^2$.
\bigskip

\noindent
\underline{\textit{Notion of  $ {\mit \Gamma} $-convergence. (cf. \cite{GammaConv})}}
Let $ X $ be a reflexive Banach space. We say that a sequence 
$ \{ \Psi_\nu  \}_{\nu>0}$ of proper functionals $ \Psi_\nu  : X \longrightarrow (-\infty, \infty] $, $ \nu > 0 $ (resp. a sequence $ \{ \Psi_n  \}_{n = 1}^\infty $ of proper functionals $ \Psi_n  : X \longrightarrow (-\infty, \infty] $, $ n \in \N $), $\Gamma$-converges to a proper functional $ \Psi_0 : X \longrightarrow (-\infty, \infty] $ as
$ \nu \downarrow 0 $
(resp. as $n\to+\infty$), if and only if the following two conditions hold:
\begin{description}
\item[\textmd{(\hypertarget{gamma1}{$\gamma\,1$})}](lower bound) $ \displaystyle \liminf_{\nu \downarrow 0} \Psi_\nu(v_\nu) \geq \Psi_0(v_0) $, if $ v_0 \in X $, $ \{ v_\nu \}_{\nu > 0} \subset X $, and $ v_\nu\to v_0 $ in $ X $ as 
$ \nu \downarrow 0 $
(resp. replacing ``$\nu$'' by ``$n$'', and ``$\nu \downarrow 0$'' by ``$n\to \infty$'');
\item[\textmd{(\hypertarget{gamma2}{$\gamma\,2$})}](optimality) for any $ v_0 \in D(\Psi_0) $, there exists a net $ \{ v_\nu \}_{\nu > 0} \subset X  $ (resp. a sequence $\{ v_n \}_{n=1}^\infty \subset X$) , such that $ v_\nu \to v_0 $ and $ \Psi_\nu(v_\nu) \to \Psi_0(v_0) $, 
$ \nu \downarrow 0 $
(resp. replacing ``$\nu$'' by ``$n$'', and ``$\nu \downarrow 0$'' by ``$n\to \infty$'').
\end{description}
\bigskip

\noindent
\underline{\textit{Notations of basic elliptic operators.}}
We denote by $\langle {}\cdot{}, {}\cdot{} \rangle_{\ast}$ the duality pairing between $H^{1}(\Omega)$ and its dual $H^{1}(\Omega)^{\ast}$. Besides, let $F\ :\ H^{1}(\Omega) \longrightarrow H^{1}(\Omega)^{\ast}$ be the duality mapping, defined as
\begin{equation}
\langle Fw, v \rangle_{\ast} := (w,v)_{H^{1}(\Omega)} = (w,v)_{L^{2}(\Omega)} + (\nabla w, \nabla v)_{L^{2}(\Omega)^{N}}, \mbox{ for all $v,w \in H^{1}(\Omega)$.}
\end{equation}
As is well-known, $Fu = - \laplace_{\rm N}u + u$ in $L^{2}(\Omega)$, if $u \in H^{1}(\Omega)$ belongs to the domain
\begin{equation*}
D_{\rm N} := \{ v \in H^{2}(\Omega)\ |\ \nabla v \cdot \nu_{\partial\Omega} = 0 \mbox{ in } L^{2}(\partial\Omega)\},
\end{equation*}
of the Laplacian operator
\begin{equation*}
\laplace_{\rm N}\ : \ u \in D_{\rm N} \subset L^{2}(\Omega) \mapsto \laplace u \in L^{2}(\Omega),
\end{equation*}
subject to the Neumann-zero boundary condition.

\bigskip
\noindent
\underline{\textit{Notations in basic measure theory. (cf. \cite{AFP})}}
 Given $ N \in \N $, we denote by $\mathcal{L}^{N}$ the $N$-dimensional Lebesgue measure, and for a measurable function $ f : B \longrightarrow [-\infty, \infty] $ on a Borel set $ B \subset \R^N $, we denote by $ [f]^+ $ and $ [f]^- $, respectively, the positive part and the negative part of $ f $. In particular, the measure theoretical phrases, such as ``a.e.'', ``$dt$'' and ``$dx$'', and so on, are all with respect to the Lebesgue measure in each corresponding dimension.

    For any open set $ U \subset \mathbb{R}^N $, we denote by $ \mathcal{M}(U) $ (resp. $ \mathcal{M}_{\rm loc}(U) $) the space of all finite Radon measures (resp. the space of all Radon measures) on $ U $. In general, the space $ \mathcal{M}(U) $ (resp. $ \mathcal{M}_{\rm loc}(U) $) is known as the dual of the Banach space $ C_0(U) $ (resp. dual of the locally convex space $ C_{\rm c}(U) $), for any open set $ U \subset \mathbb{R}^d $.
\bigskip

\noindent
\underline{\textit{Notations in BV-theory. (cf. \cite{AFP, ABM, EG, G})}}
Let $ N \in \N $ be a fixed number, and let $U\subset \R^N$ be an open set. A function $ v \in L^1(U) $ (resp. $ v \in L_{\rm loc}^1(U) $)  is called a function of bounded variation, or a BV-function, (resp. a function of locally bounded variation or a BV$\empty_{\rm loc}$-function) on $ U $, if and only if its distributional differential $ D v $ is a finite Radon measure on $ U $ (resp. a Radon measure on $ U $), namely $ D v \in \mathcal{M}(U) $ (resp. $ D v \in \mathcal{M}_{\rm loc}(U) $).
We denote by $ BV(U) $ (resp. $ BV_{\rm loc}(U) $) the space of all BV-functions (resp. all BV$\empty_{\rm loc}$-functions) on $ U $. For any $ v \in BV(U) $, the Radon measure $ D v $ is called the variation measure of $ v $, and its  total variation $ |Dv| $ is called the total variation measure of $ v $. Additionally, the value $|Dv|(U)$, for any $v \in BV(U)$, can be calculated as follows:
\begin{equation*}
|Dv|(U) = \sup \left\{ \begin{array}{l|l}
\ds \int_{U} v \ {\rm div} \,\varphi \, dy & \varphi \in C_{\rm c}^{1}(U)^N \ \ \mbox{and}\ \ |\varphi| \le 1\ \mbox{on}\ U
\end{array}
\right\}.
\end{equation*}
The space $BV(U)$ is a Banach space, endowed with the following norm:
\begin{equation*}
|v|_{BV(U)} := |v|_{L^{1}(U)} + |D v|(U),\ \ \mbox{for any}\ v\in BV(U).
\end{equation*}
We say that a sequence $\{ v_{n} \}_{n = 1}^\infty \subset BV(U)$ strictly converges in $BV(U)$ to $v \in BV(U)$ if $v_{n} \to v$ in $L^{1}(U)$ and $|Dv_{n}|(U) \to |Dv|(U)$ as $ n \to \infty $. In particular, if the boundary $\partial U$ is Lipschitz, then the space $BV(U)$ is continuously embedded into $L^{N/(N-1)}(U)$ and compactly embedded into $L^{q}(U)$ for any $1 \le q < N/(N-1)$ (cf. \cite[Corollary 3.49]{AFP} or \cite[Theorem 10.1.3-10.1.4]{ABM}). Additionally, if $1 \le r < \infty$, then the space $C^{\infty}(\overline{U})$ is dense in $BV(U) \cap L^{r}(U)$ for the intermediate convergence (cf. \cite[Definition 10.1.3. and Theorem 10.1.2]{ABM}), i.e. for any $v \in BV(U) \cap L^{r}(U)$, there exists a sequence $\{v_{n} \}_{n = 1}^\infty \subset C^{\infty}(\overline{U})$ such that $v_{n} \to v$ in $L^{r}(U)$ and $\int_{U}|\nabla v_{n}|dx \to |Dv|(U)$ as $n \to \infty$.
\bigskip

\noindent
\underline{\textit{Weighted total variation and some extensions. (cf. \cite{AB, AFP, BBF})}} In this paper, we let
\begin{equation}\label{2.2}
\left\{ \begin{array}{l}
\ds X_{\rm c}(\Omega) := \left\{ \begin{array}{l|l}
\varphi \in L^{\infty}(\Omega)^{N} & {\rm div} \, \varphi \in L^{2}(\Omega) \mbox{ and supp $\varphi$ is compact in } \Omega
\end{array} \right\}, \vspace{3mm}\\
\ds W_{0}(\Omega) := \left\{ \begin{array}{l|l}
\beta \in H^{1}(\Omega) \cap L^{\infty}(\Omega) & \beta \ge 0 \mbox{ a.e. in } \Omega
\end{array} \right\}, \vspace{3mm}\\
\ds W_{\rm c}(\Omega) := \left\{ \begin{array}{l|l}
\beta \in H^{1}(\Omega) \cap L^{\infty}(\Omega) &
\parbox{5.5cm}{there exists $ \delta_{\beta} > 0 $, such that $ \beta \ge \delta_{\beta} $  a.e. in $ \Omega $}
\end{array} \right\}.
\end{array} \right.
\end{equation}
Given $\beta \in W_{0}(\Omega)$, we define a functional $\Phi_{0}(\beta;\ \cdot\ )$ on $L^{2}(\Omega)$, by letting
\begin{equation*}
\ds v \in L^{2}(\Omega) \mapsto \Phi_{0}(\beta;v) := \mbox{sup} \left\{
\begin{array}{l|l}
\ds \int_{\Omega} v\ \mbox{div}\varphi \ dx
& \parbox{3cm}{$ \varphi \in X_{\rm c}(\Omega)$ and $|\varphi| \le \beta$ a.e. in $\Omega$}
\end{array}
\right\}.
\end{equation*}

In \cite[Theorem 5]{BBF}, it is proved that the functional $\Phi_{0}(\beta;\ \cdot\ )$ coincides with the lower semi-continuous envelope of the functional $\Phi_{0}^{\circ}(\beta;\ \cdot\ )$ defined as \begin{equation*}
v \in W^{1,1}(\Omega) \cap L^{2}(\Omega) \mapsto \Phi_{0}^{\circ}(\beta;v) := \int_{\Omega} \beta|\nabla v| \ dx,
\end{equation*}on $L^{2}(\Omega)$, and therefore, this is a maximal extension of $\Phi_{0}^{\circ}(\beta;\ \cdot\ )$ in the class of proper l.s.c. and convex functions on $L^{2}(\Omega)$. Also, for any $\beta \in W_{0}(\Omega)$ and any $v \in BV(\Omega) \cap L^{2}(\Omega)$, there exists (see \cite[Theorem 4.3]{AB} and \cite[Proposition 5.48]{AFP}) a Radon measure $|Dv|_{\beta} \in \M(\Omega)$ such that $|D v|_{\beta}$ is absolutely continuous with respect to $|Dv|$ and
\begin{equation*}
\Phi_{0}(\beta;v) = |Dv|_{\beta}(\Omega) = \int_{\Omega}d|Dv|_{\beta}.
\end{equation*}
Moreover, by \cite[Theorem 6.1]{AdCF}, for any $\beta\in W_0(\Omega)$ and any open set $U \subset \Omega$, it holds that
{
\begin{equation}\label{2.4}
\begin{array}{rcl}
|Dv|_{\beta}(U)
& = & \mbox{inf} \left \{
\begin{array}{l|l}
\ds \liminf_{n \to \infty} \int_{U} \beta |\nabla v_{n}| dx
& \parbox{5.5cm}{ $\{v_{n}\}_{n = 1}^\infty \subset W^{1,1}(U) \cap L^{2}(U)$ \mbox{ such that } $v_{n} \to v$ \mbox{ in } $L^{2}(U)$ \mbox{ as } $n \to \infty$}
\end{array}
\right\}
\\
& = & \ds \int_U \widetilde{\beta} d|Dv|,
\end{array}
\end{equation}
where $ \widetilde{\beta} $ denotes the continuous representative of $\beta$.
}
\bigskip

\noindent
\underline{\textit{Generalized weighted total variation.}}
For any $\beta \in H^{1}(\Omega) \cap L^{\infty}(\Omega)$ and any $v \in BV(\Omega) \cap L^{2}(\Omega)$, we define a real Radon measure $[\beta|Dv|] \in \M(\Omega)$, as follows:
{
\begin{equation*}
\begin{array}{rcl}
[\beta|D v|](B) & := & |D v|_{[\beta]^{+}}(B) - |Dv|_{[\beta]^{-}}(B)
\\[1ex]
& = & \ds \int_B(\widetilde{[\beta]^+}-\widetilde{[\beta]^-})d|Dv|, \mbox{ \ for any Borel set } B \subset \Omega.
\end{array}
\end{equation*}
}
Then, {$[\beta|D v|](\Omega)$} might be called the total variation of $v \in BV(\Omega) \cap L^{2}(\Omega)$ weighted by the generally sign-changing weight $\beta \in H^{1}(\Omega) \cap L^{\infty}(\Omega)$. In this paper, we call the value $[\beta|D v|](\Omega)$ ``the generalized weighted total variation" in short.
\begin{rem}\label{rem-wtv}
\begin{em}
With regard to the generalized weighted total variations, the following facts are verified in \cite{MS}:
\begin{description}
\item [{(\hypertarget{Fact1}{Fact\,1})}] (Strict approximation, {cf. \cite[Lemma 1]{MS}}) Let $\beta \in H^{1}(\Omega) \cap L^{\infty}(\Omega)$ and $v \in BV(\Omega) \cap L^{2}(\Omega)$ be arbitrarily fixed functions, and let $\{ v_{n} \}_{n = 1}^\infty \subset C^{\infty}(\overline{\Omega})$ be a sequence such that
\begin{equation*}
v_{n} \to v \mbox{ in } L^{2}(\Omega) \mbox{ and strictly in } BV(\Omega) \mbox{ as } n \to \infty.
\end{equation*}
Then,
\begin{equation*}
\int_{\Omega} \beta |\nabla v_{n}| dx \to \int_{\Omega} d[\beta|Dv|] \mbox{ as } n \to \infty.
\end{equation*}
\item [(\hypertarget{Fact2}{Fact\,2})]{(cf. \cite[Lemma 2]{MS})} For any $v\in BV(\Omega) \cap L^{2}(\Omega)$, the mapping
\begin{equation*}
\beta \in H^{1}(\Omega) \cap L^{\infty}(\Omega) \mapsto \int_{\Omega}d[\beta|Dv|] \in \R,
\end{equation*}
is a linear functional. Moreover, if $\varphi \in H^{1}(\Omega) \cap C(\overline{\Omega})$ and $\beta \in H^{1}(\Omega) \cap L^{\infty}(\Omega)$, then,
\begin{equation*}
\int_{\Omega} d[\varphi\beta|D v|] = \int_\Omega  \varphi \, d[\beta|D v|].
\end{equation*}
\end{description}
\end{em}
\end{rem}
\bigskip

Next, we consider the time-dependent cases of the weighted-total variations. Here we collect some additional notations for the convenience of descriptions.
\bigskip

\noindent
\underline{\textit{Time-dependent weighted total variation.}}
Let $I \subset (0,\infty)$ be a bounded open interval. In what follows, we set
\begin{equation*}
\begin{array}{ll}
\ds \mathscr{L}^{2}_{0}(I;\Omega) := \left\{ \begin{array}{l|l}
\beta \in L^{2}(I;L^{2}(\Omega)) & \beta(t) \ge 0, \mbox{ \ a.e. in } \Omega, \mbox{ a.e. } t \in I
\end{array} \right\}, \vspace{3mm}\\
\ds \mathscr{W}_{0}(I;\Omega) := \left\{ \begin{array}{l|l}
\beta \in L^{\infty}(I;H^{1}(\Omega)) & \beta(t) \in W_{0}(\Omega), \mbox{ a.e. } t \in I
\end{array} \right\}, \vspace{3mm}\\
\ds \mathscr{W}_{\rm c}(I;\Omega) := \left\{ \begin{array}{l|l}
\beta \in L^{\infty}(I;H^{1}(\Omega)) & \beta(t) \in W_{\rm c}(\Omega), \mbox{ \ a.e. } t \in I
\end{array} \right\}.
\end{array}
\end{equation*}
Given $\beta \in \mathscr{W}_{0}(I;\Omega)$, we define a functional $\Phi_{0}^{I}(\beta;\ \cdot\ )$ on $L^{2}(I;L^{2}(\Omega))$ by letting
\begin{equation}\label{2.4-1}
v \in L^{2}(I;L^{2}(\Omega)) \mapsto \Phi_{0}^{I}(\beta;v) := \left \{
\begin{array}{ll}
\multicolumn{2}{l}{
\ds \int_{I} \Phi_{0}(\beta(t) ; v(t) ) dt,
}
\\[1ex]
& \mbox{if } \Phi_{0}(\beta(\cdot);v(\cdot)) \in L^{1}(I),
\\[2ex]
\ds \infty, & \mbox{otherwise,}
\end{array}
\right.
\end{equation}
and we call this functional the \em time-dependent weighted total variation.\em
\medskip

\begin{rem}\label{rem-2}
\begin{em}
With regard to time-dependent weighted total variations, the following facts are verified in \cite{MS,SW2014}:
\begin{description}
\item[(\hypertarget{Fact3}{Fact\,3})]{(cf. \cite[Lemma 4]{MS})} The following three items hold:
\item[\textmd{~~~~$\bullet$}]If $\beta \in \mathscr{L}^{2}_{0}(I; \Omega) $ and $v \in L^{2}(I;H^{1}(\Omega))$, then the function $ t \in I \mapsto \Phi_{0}(\beta(t);v(t))$ is integrable on $I$,
\item[\textmd{~~~~$\bullet$}]If $\beta \in \mathscr{W}_{\rm c}(I;\Omega) \cap C(\overline{I};L^{2}(\Omega))$, $v \in C(\overline{I};L^{2}(\Omega))$, then the function $t \in I \mapsto \Phi_{0}(\beta(t);v(t))$ is l.s.c. on $I$,
\item[\textmd{~~~~$\bullet$}]If $\beta \in \mathscr{W}_{0}(I;\Omega) \cap C(\overline{I};L^{2}(\Omega))$, $v \in C(\overline{I};L^{2}(\Omega))$ and $v(t) \in BV(\Omega)$ a.e. $t \in I$, then the function $t \in I \mapsto \Phi_{0}(\beta(t);v(t))$ is measurable on $I$.
\item[(\hypertarget{Fact6}{Fact\,4})]{(Strict approximation, cf. \cite[Lemma 5]{MS}, \cite[Remark 2]{SW2014})} Given any $v \in  L^{2}(I;L^{2}(\Omega))$ such that $|D v(\cdot)|(\Omega) \in L^{1}(I)$, there exists a sequence $\{ v_{n} \}_{n = 1}^\infty \subset C^{\infty}(\overline{I\times\Omega})$ of smooth functions, such that
\begin{equation*}
v_{n} \to v \mbox{ in } L^{2}(I;L^{2}(\Omega)), ~~ \int_{I} \left| \int_{\Omega} |\nabla v_{n}(t)| dx -\int_{\Omega}d|D v(t)| \right| dt \to 0,
\end{equation*}
\vspace{-2ex}
{
\begin{equation*}
v_n(t) \to v(t) \mbox{ in } L^2(\Omega) \mbox{ \ and \ } \int_{\Omega} |\nabla v_{n}(t)| dx \to \int_{\Omega}d|D v(t)|, \mbox{ a.e. $ t \in I $, \ as $ n \to \infty $.}
\end{equation*}
}
\item[(\hypertarget{Fact7}{Fact\,5})]{Let $ \beta \in \mathscr{W}_0(I; \Omega) \cap C(\overline{I}; L^2(\Omega)) \cap L^\infty(I \times \Omega) $, $ \{ \beta_n \}_{n = 1}^\infty \subset \mathscr{L}_0^2(I; \Omega) $, $ v \in C(\overline{I}; L^2(\Omega)) $ and $ \{ v_n \}_{n = 1}^\infty \subset L^2(I; H^1(\Omega)) $  be such that
\begin{equation}\label{fact7-0}
\left\{ ~ \parbox{8.4cm}{
$ \beta_n(t) \to \beta(t) $ in $ L^2(\Omega) $ and weakly in $ H^1(\Omega) $,
\\[1ex]
$ v_n(t) \to  v(t) $ in $ L^2(\Omega) $,
} \right. \mbox{a.e. $ t \in I $, as $ n \to \infty $,}
\end{equation}
and
\begin{equation}\label{fact7-1}
\parbox{12.5cm}{
$ \beta \geq \delta_0 $ and $ \ds \inf_{n \in \N} \beta_n \geq \delta_0 $, a.e. in $ I \times \Omega $, for some constant $ \delta_0 > 0 $.
}
\end{equation}
Let $ \varrho \in C(\overline{I}; L^2(\Omega)) \cap L^\infty(I; H^1(\Omega)) \cap L^\infty(I \times \Omega) $ and $ \{ \varrho_n \}_{n = 1}^\infty \subset L^2(I; L^2(\Omega)) $ be such that
\begin{equation}\label{fact7-2}
\left\{ \parbox{12.5cm}{
$ \varrho_n(t) \to \varrho(t) ${ in $ L^2(\Omega) $ and weakly in $ H^1(\Omega) $,  a.e. $ t \in I $, as $ n \to \infty $,}
\\[1ex]
$ |\varrho| \leq M_0 $ and $ \ds \sup_{n \in \N} |\varrho_n| \leq M_0 $, a.e. in $ I \times \Omega $, for some constant $ M_0 > 0 $.
} \right.
\end{equation}
In addition, let us assume that
\begin{equation*}
\ds \int_I \int_\Omega \beta_n(t) |\nabla v_n(t)| \, dx dt \to \int_I \int_\Omega d \bigl[ \beta(t) |Dv(t)| \bigr] \, dt \mbox{ \ as $ n \to \infty $.}
\end{equation*}
Then,
}
\begin{equation*}
\int_I \int_\Omega \varrho_n(t) |\nabla v_n(t)| \, dx dt \to \int_I \int_\Omega d \bigl[ \varrho(t) |D v(t)| \bigr] \, dt \mbox{ \ as $ n \to \infty $.}
\end{equation*}
\end{description}
\end{em}
\end{rem}


\hypertarget{Sec.2}{\empty}
\section{Statement of the main results}
\ \ \vspace{-3ex}

Here we list the assumptions for the system (\hyperlink{(S)}{S}):
\begin{enumerate}
\item[\hypertarget{(H1)}{(H1)}]$g : \R \longrightarrow \R$ is a locally Lipschitz continuous function such that $ g(0) \leq 0 $,  $ g(1) \geq 0 $, and $ g $ has a nonnegative primitive $\hat{g} \in W^{2,\infty}_{\rm loc}(\R)$.
\item[\hypertarget{(H2)}{(H2)}]$\alpha_0 : \R \longrightarrow (0, \infty) $ is a locally Lipschitz continuous function.
\item[\hypertarget{(H3)}{(H3)}]$\alpha : \R \longrightarrow (0, \infty) $ is a $ C^{2} $-function, such that  $ \alpha $ is convex on $ \R $, and $ \alpha'(0) = 0 $.
\item[\hypertarget{(H4)}{(H4)}] There exists a positive constant $\delta_{\alpha} \in (0, 1) $, such that
\begin{equation*}
\alpha_{0}(\tau) \ge \delta_{\alpha} \mbox{ and } \alpha(\tau) \ge \delta_{\alpha}, \mbox{ for all } \tau \in \R.
\end{equation*}
\item[\hypertarget{(H5)}{(H5)}] The initial data $ [\eta_0, \theta_0] $ belongs to a class $ D_0 \subset L^2(\Omega)^2 $, defined as:
\begin{equation*}
D_0 := \left\{ \begin{array}{l|l}
[w, v] \in L^2(\Omega)^2 & 0 \leq w \leq 1 \mbox{ a.e. in $ \Omega $ \ and \ $v \in L^{\infty}(\Omega)$}
\end{array} \right\}.
\end{equation*}
\end{enumerate}
\begin{rem}[Possible choice of given functions]
\begin{em}
Referring to \cite{KWC1, KWC2}, the setting
\begin{equation*}
\ds g(\tau) = \tau -1 \mbox{ with } \hat{g}(\tau) := \frac{1}{2}(\tau-1)^{2} \mbox{ and } \alpha_{0}(\tau) = \alpha(\tau) = \tau^{2} + \delta_{\alpha}, \mbox{ for } \tau \in \R,
\end{equation*}
provides a possible choice of given functions that fulfills the above assumptions.
\end{em}
\end{rem}
\begin{rem}\label{Rem.FEreal}
\begin{em}
On the basis of the assumptions (H1)-(H4) and the notations prepared in Section 1, the exact formulation of the free-energy $\mathscr{F}$ in (\ref{free}) can be prescribed as follows:
\begin{equation}\label{freeEnergyReal}
[\eta, \theta] \in L^{2}(\Omega) \times L^{2}(\Omega) \mapsto \mathscr{F}(\eta,\theta) := \left\{
\begin{array}{ll}
\multicolumn{2}{l}{
\ds \frac{1}{2} \int_{\Omega} |\nabla\eta|^{2} dx + \int_{\Omega} \hat{g}(\eta) dx + \Phi_{0}(\alpha(\eta);\theta),
}
\\[3ex]
& \ds \mbox{if } \eta \in H^{1}(\Omega) \cap L^{\infty}(\Omega)
\\[1ex]
& \ds \mbox{and } \theta \in D(\Phi_{0}(\alpha(\eta);\ \cdot\ )),
\\[2ex]
\ds \infty,  & \mbox{otherwise. }
\end{array}
\right.
\end{equation}
\end{em}
\end{rem}

First of all, we give the exact definition of the solution to the system (\hyperlink{(S)}{S}).
\begin{defn}[Definition of solution to (\hyperlink{(S)}{S})]\label{sol}
\begin{em}
A pair $[\eta,\theta]$ of functions $ \eta = $\linebreak $ \eta(t, x) $ and $ \theta =  \theta(t, x) $ is called an energy-dissipative solution to (\hyperlink{(S)}{S}), if and only if the components $\eta$ and $\theta$ fulfill the following four conditions.
\begin{enumerate}
\item[\hypertarget{(S0)}{(S0)}] $\eta \in C([0,\infty);L^{2}(\Omega)) \cap W^{1,2}_{\rm loc}((0,\infty);L^{2}(\Omega)) \cap L_{\rm loc}^2([0, \infty); H^1(\Omega)) \cap L^{\infty}_{\rm loc}((0,\infty);H^{1}(\Omega)) $ $ \cap L^\infty(Q) $, and 
$ 0 \le \eta \le 1 $  a.e. in  $Q$;
\\[1ex]
$\theta \in C([0,\infty);L^{2}(\Omega)) \cap W^{1,2}_{\rm loc}((0,\infty);L^{2}(\Omega)) \cap L^{\infty}(Q) $, $ |D \theta({}\cdot{})|(\Omega) \in L_{\rm loc}^1([0, \infty)) \cap L_{\rm loc}^{\infty}((0,\infty))$, and
$ |\theta| \le |\theta_{0}|_{L^{\infty}(\Omega)} $ a.e. in $ Q $.
\item[\hypertarget{(S1)}{(S1)}] $\eta$ solves the following variational identity:
\begin{equation}\label{vari01}
\hspace{-2ex}
\begin{array}{c}
\ds \int_{\Omega} \left( \rule{0pt}{10pt} \eta_{t}(t) + g(\eta(t)) \right) w \, dx + \int_{\Omega}\nabla\eta(t) \cdot \nabla w \, dx + \int_{\Omega}d[w \alpha'(\eta(t))|D \theta(t)|] = 0,
\\[2ex]
\mbox{for any $w \in H^{1}(\Omega) \cap L^{\infty}(\Omega)$ and a.e. $t \in (0,\infty)$,}
\end{array}
\end{equation}
subject to the initial condition \ $ \eta(0) = \eta_0 $  in $ L^2(\Omega) $.
\item[\hypertarget{(S2)}{(S2)}]
$\theta$ solves the following variational inequality:
\begin{equation}\label{vari02}
\begin{array}{c}
\ds \int_{\Omega} \alpha_{0}(\eta(t))\theta_{t}(t) (\theta(t)-v) \, dx + \Phi_{0}(\alpha(\eta(t));\theta(t)) \le \Phi_{0}(\alpha(\eta(t));v),
\\[2ex]
\mbox{for any $v \in BV(\Omega) \cap L^{2}(\Omega)$ and a.e. $t \in (0,\infty)$,}
\end{array}
\end{equation}
subject to the initial condition \
$ \theta(0) = \theta_0 $ in $ L^2(\Omega) $.
\item[\hypertarget{(S3)}{(S3)}]
There exists a function $\mathscr{J}_* \in L^{1}_{\rm loc}([0,\infty)) \cap BV_{\rm loc}((0,\infty))$, such that $\mathscr{J}_*$ is nonincreasing on $(0,\infty)$, and $ \mathscr{J}_*(t) = \mathscr{F}(\eta(t),\theta(t))$ for a.e. $ t \in (0, \infty) $.
\end{enumerate}
\end{em}
\end{defn}
\begin{rem}\label{rem-ken100}
\begin{em}
By (\hyperlink{Fact2}{Fact\,2}), the variational identity (\ref{vari01}) in Definition \ref{sol}, can be rewritten into the following weak formulation:
$$
\eta_{t}(t) + F\eta(t) - \eta(t) + g(\eta(t)) + [\alpha'(\eta(t))|D \theta(t)|] = 0
\mbox{ \ in}\ H^{s}(\Omega)^{\ast},\ \ \mbox{a.e.}\ t\in(0,\infty),
$$
where $s > N/2$ is a large constant such that the embedding $ H^s(\Omega) \subset C(\overline{\Omega}) $ holds true.
So, in accordance with Definition \ref{sol}, the nonstandard term $\alpha'(\eta)|D\theta|$ as  in (\ref{1steq}) is equated with the element $[\alpha'(\eta(\cdot))|D \theta(\cdot)|]  \in L^{2}_{\rm loc}((0,\infty) ; H^{s}(\Omega)^{\ast})$, and moreover, $[\alpha'(\eta(t))|D\theta(t)|] \in H^{s}(\Omega)^{\ast} \cap \M(\Omega)$ for a.e. $ t \in (0,\infty)$.

\medskip
Meanwhile, $\theta$ solves the following inclusion:
$$
\alpha_{0}(\eta(t))\theta_{t}(t) + \partial\Phi_{0}(\alpha(\eta(t));\theta(t)) \ni 0 \ \ \mbox{in}\ L^{2}(\Omega), \ \ \mbox{a.e.}\ t\in(0,\infty).
$$
Therefore, the mathematical meaning of the nonstandard term $-\mbox{div}(\alpha(\eta)\frac{D\theta}{|D\theta|})$ in $(\ref{vari01})$ is given in terms of the subdifferential $ \partial \Phi_{0}(\alpha(\eta(t));{}\cdot{})$ of the unknown-dependent total variation.
\end{em}
\end{rem}
\medskip

On the basis of Definition \ref{sol}, the main results in this paper are stated as follows.
\begin{mainTh}[Solvability of (\hyperlink{(S)}{S})]\label{mainTh1}
Let us assume (\hyperlink{(H1)}{H1})-(\hyperlink{(H4)}{H4}). Then, there exists at least one energy-dissipative solution $[\eta,\theta]$ to the system (\hyperlink{(S)}{S}).

\end{mainTh}
\begin{mainTh}[Large-time behavior]\label{mainTh2}
Under assumptions (\hyperlink{(H1)}{H1})-(\hyperlink{(H4)}{H4}), let $[\eta,\theta]$ be an energy-dissipative solution to (\hyperlink{(S)}{S}), and let $ \omega(\eta, \theta) $ be the $ \omega $-limit set of $ [\eta, \theta] $, i.e.:
\begin{equation*}
\omega(\eta, \theta) := \left\{ \begin{array}{l|l}
[\eta_\infty, \theta_\infty] \in L^2(\Omega)^2 & \parbox{6.4cm}{
$ [\eta(t_n), \theta(t_n)] \to [\eta_\infty, \theta_\infty] $ in $ L^2(\Omega)^2 $ as $ n \to \infty $, \ for some $ \{ t_n \}_{n = 1}^\infty \subset (0, \infty) $ satisfying $ t_n \uparrow \infty $ as $ n \to \infty $
}
\end{array} \right\}.
\end{equation*}
Then, the following two items hold.
\begin{enumerate}
\item[(\hypertarget{(O)}{O})]$ \omega(\eta, \theta) \ne \emptyset $, and $ \omega(\eta, \theta) $ is compact in $ L^2(\Omega)^2 $.
\item[(\,\hypertarget{(I)}{I}\,)]Any  $ \omega $-limit point, $ [\eta_\infty, \theta_\infty] \in \omega(\eta, \theta) $, fulfills that:
\begin{enumerate}
\item[(\hypertarget{i-a}{i-a})]$ 0 \leq \eta_\infty \leq 1 $ and $ |\theta_\infty| \leq |\theta_0|_{L^\infty(\Omega)} $ a.e. in $ \Omega $;
\item[(\hypertarget{i-b}{i-b})]$ -{\mit \Delta}_{\rm N} \eta_\infty +g(\eta_\infty) = 0 $ in $ L^2(\Omega) $;
\item[(\hypertarget{i-c}{i-c})]$ \theta_\infty $ is a constant over $ \Omega $, i.e. $ \theta_\infty $ is a global minimizer of the convex function $ \Phi_0(\alpha(\eta_\infty);{}\cdot{}) $ on $ L^2(\Omega) $.
\end{enumerate}
\end{enumerate}
\end{mainTh}

\section{Approximating problems}
\ \ \vspace{-3ex}

In this Section, we prescribe the approximating problems to the system (\hyperlink{S}{S}), and verify some  key-properties of the approximating solutions.
To this end, we take a ``{suitable approximation} $ \{ |{}\cdot{}|_\nu \}_{\nu \in (0, 1)} \subset C^1(\R) $'' of the Euclidean norm, and for any  $ 0 \leq \beta \in L^2(\Omega) $, we define $ \{ \Phi_{\nu}(\beta;{}\cdot{}) \}_{\nu \in (0, 1)} $ a collection of proper l.s.c. and convex functions on $ L^2(\Omega) $, by putting
\begin{equation}\label{relaxPhi}
z \in L^2(\Omega) \mapsto \Phi_{\nu}(\beta; v) := \left\{ \begin{array}{ll}
\multicolumn{2}{l}{\ds \int_\Omega  \beta|\nabla v|_\nu \, dx +\frac{\nu}{2} \int_\Omega |\nabla v|^2 \, dx,}
\\[2ex]
& \mbox{if $ v \in H^1(\Omega) $,}
\\[2ex]
\infty, & \mbox{otherwise,}
\end{array} \right. \mbox{for any $ \nu \in (0, 1) $.}
\end{equation}
Observe that the convex function $ \Phi_{\nu}(\beta;{}\cdot{}) $ corresponds to a relaxed version of the weighted-total variation $ \Phi_0(\beta;{}\cdot{}) $. The precise definition of a suitable approximation is the following one:

\begin{defn}\label{euclregul}
\begin{em}
We say that a collection of functions $\{|\cdot|_\nu\}_{\nu\in(0,1)}$ is a suitable approximation to the {Euclidean} norm if the following properties hold:
\begin{description}
\item [\textmd{(\hypertarget{AP1}{AP1})}] $|\cdot|_\nu:\R^N\mapsto [0,+\infty[$ is a convex $C^1$ function such that $|0|_\nu = 0 $, for all $ \nu \in (0, 1) $.

\item [\textmd{(\hypertarget{AP2}{AP2})}] There exist bounded functions $ q_{0}\ :\ (0,1) \longrightarrow (0,1]$, $ q_{1} : (0, 1) \longrightarrow [1, \infty) $, $ r_{k}\ :\ (0,1) \longrightarrow [0,\infty)$, $k = 0, 1$, such that:
\begin{equation*}
q_0(\nu) \to 1, ~ q_1(\nu) \to 1, ~ r_0(\nu) \to 0 \mbox{ and } r_1(\nu) \to 0, \mbox{ as $ \nu \downarrow 0 $,}
\end{equation*}
and
\begin{equation*}
\begin{array}{c}
\ds
|\xi|_\nu \geq q_0(\nu) |\xi| -r_0(\nu) \mbox{ and } |[\nabla |\cdot|_\nu](\xi)| \leq q_1(\nu) |\xi|^{r_1(\nu)}
\\[1ex]
\mbox{for any $ \xi \in \R^N $ and $ \nu \in (0, 1) $.}
\end{array}
\end{equation*}
\end{description}
\end{em}
\end{defn}

{
\begin{rem}\label{Rem.Euc}
\begin{em}
Note that (\hyperlink{AP1}{AP1})-(\hyperlink{AP2}{AP2}) lead to the following fact:
\begin{center}
$ |\xi|_\nu \leq [\nabla |\cdot|_\nu](\xi) \cdot \xi \leq q_1(\nu)|\xi|^{1 +r_1(\nu)} $, \ for all $ \xi \in \R^N $ and $ \nu \in (0, 1) $. \end{center}
Also, we note that the class of possible regularizations verifying  (\hyperlink{AP1}{AP1})-(\hyperlink{AP2}{AP2}) covers a number of standard type regularizations. For instance,
\begin{itemize}
\item Hyperbola type, i.e. $\xi \in \R^N \mapsto \sqrt{|\xi|^2 +\nu^2} -\nu$, for $ \nu \in (0, 1) $,
\item Yosida's regularization, i.e.
$ \xi \in \R^N \mapsto |\xi|_\nu := \ds \inf_{\varsigma \in \R^N} \left\{ |\xi| +\frac{\nu}{2} |\varsigma -\xi|^2 \right\} $, for $ \nu \in (0, 1) $,
\vspace{-1ex}
\item Hyperbolic-tangent type, i.e.
$ \xi \in \R^N \mapsto |\xi|_\nu := \ds \int_0^{|\xi|} \tanh \frac{\tau}{\nu} \, d \tau $, for $ \nu \in (0, 1) $,
\vspace{-1ex}
\item Arctangent type, i.e.
$ \xi \in \R^N \mapsto |\xi|_\nu := \ds \frac{2}{\pi} \int_0^{|\xi|} {\rm Tan}^{-1} \frac{\tau}{\nu} \, d \tau $, for $ \nu \in (0, 1) $,
\vspace{-1ex}
\item $ p $-growth type, i.e.
$ \xi \in \R^N \mapsto |\xi|_\nu := \ds \frac{1}{p(\nu)} |\xi|^{p(\nu)} $, for $ \nu \in (0, 1) $,
with a function $ p : (0, 1) \longrightarrow (1, \infty) $ satisfying $ p(\nu) \downarrow 1 $ as $  \nu \downarrow 0 $.
\end{itemize}
\end{em}
\end{rem}
}

{
\begin{rem}\label{Rem.|.|_nu}
\begin{em}
Roughly summarized, the term ``suitable approximation'' means that the sequence of relaxed functionals in \eqref{relaxPhi}  $ \Gamma $-converges to the generalized weigthed total variation (see Theorems \ref{lem-gamma0} and \ref{lem-gamma}). In the proof of Main Theorems, the detailed expressions of $ \{ |\cdot|_\nu \}_{\nu \in (0, 1)} $ will not be essential, and the importance will be just in the suitability of $ \{ |\cdot|_\nu \}_{\nu \in (0, 1)} $.
\end{em}
\end{rem}
}

Now, for any $ \nu \in (0, 1) $, we can define a relaxed free-energy $ \mathscr{F}_{\nu}(\eta, \theta) $, by letting
\begin{equation}\label{relaxFreeEnrg}
[\eta, \theta] \in L^2(\Omega)^2 \mapsto \mathscr{F}_{\nu}(\eta, \theta) := \left\{ \begin{array}{ll}
\multicolumn{2}{l}{\ds \frac{1}{2} \int_\Omega |\nabla \eta|^2 \, dx +\int_\Omega \hat{g}(\eta) \, dx +\Phi_{\nu}(\alpha(\eta); \theta),}
\\[2ex]
& \mbox{if $ [\eta, \theta] \in H^1(\Omega)^2 $,}
\\[2ex]
\infty, & \mbox{otherwise.}
\end{array} \right.
\end{equation}
As a gradient system for this free-energy, we  derive the following relaxed system:
\bigskip

\noindent
$(\mbox{S$_\nu $})$:
\begin{equation*}
\ds (\eta_{\nu})_{t}(t) -{\laplace}_{\rm N}\eta_{\nu}(t) + g(\eta_{\nu}(t)) + \alpha'(\eta_{\nu}(t))|\nabla\theta_{\nu}(t)|_\nu = 0 
\ds\ \mbox{ in}\ L^{2}(\Omega),\ \mbox{ a.e. } t\in(0,\infty),
\vspace{1mm}
\end{equation*}
\begin{equation*}
\ds \alpha_{0}(\eta_{\nu}(t))(\theta_{\nu})_{t}(t) + \partial\Phi_{\nu}(\alpha(\eta_{\nu}(t));\theta_{\nu}(t)) \ni 0 \mbox{ \ in \ } L^{2}(\Omega),\ \mbox{ a.e. } t\in(0,\infty),
\vspace{1mm}
\end{equation*}
\begin{equation*}
\ds [\eta_{\nu}(0), \theta_{\nu}(0)] = [\eta_{0}^\nu, \theta_{0}^\nu] \ \mbox{ in}\ L^{2}(\Omega).
\vspace{1mm}
\end{equation*}
The initial data $[\eta_{0}^{\nu}, \theta_{0}^{\nu}]$ satisfy
\begin{equation}\label{initial-2}
[\eta_{0}^\nu, \theta_{0}^\nu] \in D_1 := D_0 \cap H^1(\Omega)^2.
\end{equation}
{Also, we need to consider} the following subset  of $ D_1 $:
\begin{equation}\label{D*(theta0)}
D_*(\theta_0) := \left\{ \begin{array}{l|l}
[w, v] \in H^1(\Omega)^2 & \parbox{6cm}{$ 0 \leq w \leq 1 $ and $ |v| \leq |\theta_0|_{L^\infty(\Omega)} $, a.e. in $ \Omega $}
\end{array} \right\}.
\end{equation}
{Additionally}, for any $\nu \in (0,1)$ and any $\beta \in L^2(\Omega)$, $\partial \Phi_{\nu}(\beta;\ \cdot\ )$ denotes the $ L^2 $-subdifferential of the convex function $\Phi_{\nu}(\beta;\ \cdot\ )$ on $L^{2}(\Omega)$.

\bigskip
For a fixed $h\in (0,1) $, we denote by $ \mbox{(AP$_h^\nu $)} $ the following problem corresponding to a particular time-discretization system for $ \mbox{(S$_\nu $)} $:
\begin{description}
\item[\textmd{$ \mbox{(AP$_h^\nu $)} $:}]for any $ [\eta_0^\nu, \theta_0^\nu] \in D_1 $, find a sequence $ \{ [\eta_{h, i}^\nu, \theta_{h, i}^\nu] \}_{i = 1}^\infty \subset H^1(\Omega)^2 $, such that
\end{description}
\vspace{-2ex}
\begin{equation*}
0 \leq \eta_{h, i}^\nu \leq 1 \mbox{ and } |\theta_{h, i}^\nu| \leq |\theta_0^\nu|_{L^\infty(\Omega)} \mbox{ a.e. in $ \Omega $, for any $ i \in \N $,}
\end{equation*}
\begin{equation}\label{kenApp01}
\begin{array}{ll}
\multicolumn{2}{l}{\ds
\frac{1}{h} (\eta_{h, i}^\nu -\eta_{h, i -1}^\nu, w)_{L^2(\Omega)} +(\nabla \eta_{h, i}^\nu, \nabla w)_{L^2(\Omega)^N} +(g(\eta_{h, i}^\nu), w)_{L^2(\Omega)}
}
\\[2ex]
\qquad \qquad & \ds +(\alpha'(\eta_{h, i}^\nu) |\nabla \theta_{h, i-1}^\nu|_\nu, w)_{L^2(\Omega)} = 0, \mbox{ for any $ w \in H^1(\Omega) $ and any $ i \in \N $,}
\end{array}
\end{equation}
\begin{equation}\label{kenApp02}
\begin{array}{c}
\ds \frac{1}{h} (\alpha_0(\eta_{h, i}^\nu)(\theta_{h, i}^\nu -\theta_{h, i -1}^\nu), v)_{L^2(\Omega)} +(\alpha(\eta_{h,i}^{\nu}) [\nabla |\cdot|_{\nu}](\nabla \theta_{h,i}^{\nu}), \nabla v)_{L^{2}(\Omega)} \vspace{3mm}
\\[2ex]
+\nu (\nabla \theta_{h, i}^\nu, \nabla v)_{L^2(\Omega)^N} = 0, \mbox{ \ for any $ v \in H^1(\Omega) $ and any $ i \in \N $,}
\end{array}
\vspace{1ex}
\end{equation}
\indent
subject to
\begin{equation*}
[\eta_{h,0}^{\nu}, \theta_{h,0}^{\nu}] = [\eta_{0}^\nu, \theta_{0}^\nu]\ \mbox{in}\ L^{2}(\Omega)^2.
\vspace{1ex}
\end{equation*}

We call the above sequence $ \{ [\eta_{h, i}^\nu, \theta_{h, i}^\nu] \}_{i = 1}^\infty \subset H^1(\Omega)^2 $ a solution to the approximating problem $ \mbox{(AP$_h^\nu $)} $, or an approximating solution in short.


\begin{thm}[{\boldmath Key-properties of $ (\mbox{AP$_h^\nu $ })$}]\label{th1}

Under assumptions (\hyperlink{(H1)}{H1})-(\hyperlink{(H4)}{H4}), there exists $ h_* \in (0, 1) $ such that, for any $ h \in (0, h_*) $ and $[\eta_0^\nu,\theta_0^\nu]\in D_1$, $(AP)_{h}^{\nu}$ admits a unique solution $ \{ [\eta_{h,i}^{\nu}, \theta_{h,i}^{\nu}]\}_{i = 1}^\infty $,  verifying
\begin{equation}\label{ene-inq}
\begin{array}{c}
\ds \frac{1}{2h}|\eta_{h,i}^{\nu} - \eta_{h,i-1}^{\nu}|_{L^{2}(\Omega)}^{2} + \frac{1}{h} \left| {\textstyle \sqrt{\alpha_{0}(\eta_{h,i}^{\nu})}}(\theta_{h,i}^{\nu} - \theta_{h,i-1}^{\nu}) \right|_{L^{2}(\Omega)}^{2} + \mathscr{F}_{\nu}(\eta_{h,i}^{\nu},\theta_{h,i}^{\nu})
\\[2ex]
\ds  \le \mathscr{F}_{\nu}(\eta_{h,i-1}^{\nu},\theta_{h,i-1}^{\nu}), \mbox{ \ for any $ i \in \N $,}
\end{array}
\end{equation}
and
\begin{equation}\label{lem2-inq}
\begin{array}{l}
\hspace{-2ex} \ds \frac{1}{2} \sum_{i = 1}^m i |\eta_{h, i}^\nu -\eta_{h, i -1}^\nu|^{2}_{L^{2}(\Omega)} + \sum_{i = 1}^m i \left| {\textstyle \sqrt{\alpha_{0}(\eta_{h, i}^\nu)}}(\theta_{h, i}^\nu -\theta_{h, i -1}^\nu) \right|_{L^{2}(\Omega)}^{2}
\\[2ex]
\hspace{6ex} \ds + mh \mathscr{F}_{\nu}(\eta_{h, i}^\nu, \theta_{h, i}^\nu) \le h \sum_{i = 1}^m \mathscr{F}_{\nu}(\eta_{h, i -1}^\nu, \theta_{h, i -1}^\nu), \mbox{ for any $ m \in \N $.}
\end{array}
\end{equation}

\end{thm}

\noindent
\textbf{Proof.}
We omit the proofs of the existence and uniqueness of approximating solutions and \eqref{ene-inq}, because those are obtained just as in \cite[Theorem 1]{MS} with slight modifications.

For \eqref{lem2-inq}, we  multiply both sides of (\ref{ene-inq}) by $ih$. Then, we see that:
\begin{equation}\label{kenApp21}
\begin{array}{ll}
\multicolumn{2}{l}{\ds \frac{i}{2} |\eta_{h, i}^\nu -\eta_{h, i -1}^\nu|_{L^2(\Omega)}^2 +i \left| \sqrt{\alpha_0(\eta_{h, i}^\nu)}(\theta_{h, i}^\nu -\theta_{h, i -1}^\nu) \right|_{L^2(\Omega)}^2}
\\[2ex]
& \quad \ds +h \bigl( i \mathscr{F}_\nu(\eta_{h, i}^\nu, \theta_{h, i}^\nu) -(i -1) \mathscr{F}_\nu(\eta_{h, i -1}^\nu, \theta_{h, i -1}^\nu) \bigr)
\\[2ex]
\leq & h \mathscr{F}_\nu(\eta_{h, i -1}^\nu, \theta_{h, i -1}^\nu), \mbox{ for any $ i \in \N $.}
\end{array}
\end{equation}
The inequality (\ref{lem2-inq}) is obtained as a summation from $1$ to $m\in \N$ in (\ref{kenApp21}). \hfill $\Box$

\begin{lem}\label{lem1}
{
Let $ \nu \in (0, 1) $, and let $ h_* $ be the constant obtained in Theorem \ref{th1}, let $ h \in (0, h_*) $ be an arbitrary time-step, and let $ \{ \eta_{h, i}^\nu, \theta_{h, i}^\nu \}_{i = 1}^\infty $ be the solution to $ \mbox{(AP$_h^\nu $)} $ with  initial data $ [\eta_0^\nu, \theta_0^\nu] \in D_*(\theta_0) $. Under these assumptions, there exist $\nu_{\ast} \in (0,1)$ and positive constants $A_{\ast}$, $B_{\ast}$, $C_{\ast} $, depending only on $ \Omega $, $ \alpha_0 $, $ \alpha $, $ \hat{g} $ and $ \theta_0 $, such that if $ h \in (0, h_*) $ and $ \nu \in (0, \nu_*) $, then the approximating solution $ \{ \eta_{h, i}^\nu, \theta_{h, i}^\nu \}_{i = 1}^\infty $ satisfies the following energy inequality:
}
\begin{equation}\label{lem1-inq}
\hspace{-1ex}
\begin{array}{l}
\hspace{-4ex} \ds \frac{1}{2}(|\eta_{h, m}^\nu -w_{0}|_{L^{2}(\Omega)}^{2} + A_{\ast} |\eta_{h, m}^\nu -v_{0}|_{L^{2}(\Omega)}^{2})
+{\frac{B_* h}{2}} \sum_{i = 1}^{m} \mathscr{F}_\nu(\eta_{h, i -1}^\nu, \theta_{h, i -1}^\nu)
\\[3ex]
\le \ds \frac{1}{2}(|\eta_{0}^\nu-w_{0}|_{L^{2}(\Omega)}^{2} + A_{\ast} |\theta_{0}^\nu-v_{0}|_{L^{2}(\Omega)}^{2}) +\frac{h}{B_*} \mathscr{F}_\nu(\eta_0^\nu, \theta_0^\nu)
\\[3ex]
\hspace{4ex} +m h C_*(1 +|w_0|_{H^1(\Omega)}^2 +|v_0|_{H^1(\Omega)}^2),
\\[3ex]
\hspace{4ex} \mbox{for any $ m \in \N $ and any $ [w_0, v_0]\in D_*(\theta_0)$.}
\end{array}
\end{equation}
\end{lem}
\bigskip

\noindent
\textbf{Proof.}
We fix the index $ i \in \N $, arbitrary, and define a large constant $ R_* > 0 $ as

\begin{equation*}
R_* := \left[ \frac{(1 +|\alpha_0|_{W^{1, \infty}(0, 1)})(1 +|\alpha|_{C^1([0, 1])})(1 +|\hat{g}|_{W^{1, \infty}(0, 1)})(1 +|\theta_0|_{L^\infty(\Omega)})(1 +\mathscr{L}^N(\Omega))}{\delta_\alpha^2} \right]^2.
\end{equation*}

We take $ w:=\eta_{h, i}^\nu -w_0 \in H^1(\Omega) $ as the test function in (\ref{kenApp01}). Then, by (\hyperlink{(H3)}{H3})-(\hyperlink{(H4)}{H4}), {H\"older} and Young's inequalities, we have
\begin{equation*}
\begin{array}{rl}
\multicolumn{2}{l}{\ds \frac{1}{2h} (|\eta_{h, i}^\nu -w_0|_{L^2(\Omega)}^2 -|\eta_{h, i -1}^\nu -w_0|_{L^2(\Omega)}^2) +\frac{1}{2}|\nabla \eta_{h, i}^\nu|_{L^2(\Omega)^N}^2}
\\[3ex]
\qquad \qquad & \ds +\int_\Omega \bigl( \alpha(\eta_{h, i}^\nu) -\alpha(w_0) \bigr) |\nabla \theta_{h, i -1}^\nu|_\nu \, dx
\leq  \frac{1}{2}|\nabla w_0|_{L^2(\Omega)^N}^2 +R_*.
\end{array}
\end{equation*}
Then,  by (\hyperlink{(H4)}{H4}) and (\hyperlink{S0}{S0}) in Definition \ref{sol},
\begin{equation}\label{kenApp17}
\begin{array}{ll}
\multicolumn{2}{l}{\ds \frac{1}{2h} (|\eta_{h, i}^\nu -w_0|_{L^2(\Omega)}^2 -|\eta_{h, i -1}^\nu -w_0|_{L^2(\Omega)}^2) +\frac{1}{2} |\nabla \eta_{h, i}^\nu|_{L^2(\Omega)^N}^2}
\\[3ex]
& \ds +\frac{\delta_\alpha}{|\alpha|_{C([0, 1])}} \int_\Omega \alpha(\eta_{h, i -1}^\nu)|\nabla \theta_{h, i -1}^\nu|_\nu \, dx
\\[5ex]
& -|\alpha|_{C([0, 1])} \bigl| |\nabla \theta_{h, i -1}^\nu|_\nu \bigr|_{L^1(\Omega)} \leq R_*(1 +|w_0|_{H^1(\Omega)}^2).
\end{array}
\end{equation}
Similarly, we take $v:=(\theta_{h, i}^\nu -v_0)/\alpha_0(\eta_{h, i}^\nu) \in H^1(\Omega) $ as the test function in (\ref{kenApp02}). Then, applying (\hyperlink{(H2)}{H2})-(\hyperlink{(H4)}{H4}), we obtain
\begin{equation}\label{kenApp13}
\begin{array}{l}
\hspace{-2ex} \ds
\frac{1}{2h}(|\theta_{h, i}^\nu -v_0|_{L^2(\Omega)}^2 -|\theta_{h, i -1}^\nu -v_0|_{L^2(\Omega)}^2)
\\[3ex]
\ds +  \frac{\delta_\alpha}{|\alpha_0|_{C([0,1])}}\int_\Omega |\nabla \theta_{h, i}^\nu|_\nu \, dx  +\nu\int_\Omega \frac{1}{\alpha_0(\eta_{h,i}^\nu)}\nabla \theta_{h,i}^\nu\cdot\nabla (\theta_{h,i}^\nu-v_0)\, dx
\\[3ex]
\ds \leq \int_\Omega \frac{{\alpha(\eta_{h,i}^\nu)}(\theta_{h,i}^\nu-v_0)\alpha_0'(\eta_{h,i}^\nu)}{\alpha_0^2(\eta_{h,i}^\nu)} [\nabla |\cdot|_{\nu}](\nabla \theta_{h,i}^{\nu})\cdot\nabla \eta_{h,i}^\nu\, dx
\\[3ex]
\ds+\nu\int_\Omega \frac{(\theta_{h,i}^\nu-v_0)\alpha_0'(\eta_{h,i}^\nu)}{\alpha_0^2(\eta_{h,i}^\nu)} \nabla \eta_{h,i}^\nu\cdot\nabla \theta_{h,i}^\nu\, dx +\frac{|\alpha|_{C([0,1])}}{\delta_\alpha}\int_\Omega  |\nabla v_0|_\nu \, dx.
\end{array}
\end{equation}
We distinguish now two cases for the different types of suitable approximations. In case that $r_1(\nu)=0$ as in (\hyperlink{AP2}{AP2}) holds (for instance, hyperbola type or Yosida's regularization), we multiply \eqref{kenApp13} by
\begin{equation}\label{kenApp14}
A_* := \frac{|\alpha_0|_{C([0, 1])}|\alpha|_{C([0, 1])}}{\delta_\alpha} ~ (\leq  \delta_\alpha R_*^{1/2}),
\end{equation} and apply Young's inequality twice with $
\ds \varepsilon_1 :=  \frac{1}{4 A_*} \mbox{ \ and \ } \ds \varepsilon_2 := \frac{1}{2|\alpha_0|_{C([0, 1])}},
$
respectively. Then, we obtain
\begin{equation*}
\begin{array}{rl}
\multicolumn{2}{l}{
\ds \frac{A_*}{2h} (|\theta_{h, i}^\nu -v_0|_{L^2(\Omega)}^2 -|\theta_{h, i -1}^\nu -v_0|_{L^2(\Omega)}^2)
}
\\[3ex]
\qquad \qquad & \ds +|\alpha|_{C([0, 1])} \bigl| |\nabla \theta_{h, i}^\nu|_\nu \bigr|_{L^1(\Omega)} +\frac{|\alpha|_{C([0, 1])}}{2 \delta_\alpha} \cdot \nu |\nabla \theta_{h, i}^\nu|_{L^2(\Omega)^N}^2
\\[3ex]\leq
& \ds  
\frac{A_* |\alpha_0|_{C[0, 1]}}{\delta_\alpha^2}
 \cdot \nu |\nabla v_0|_{L^2(\Omega)^N}^2+\left( \frac{1}{8} +4 \nu |\alpha_0|_{C([0, 1])} A_* R_* \right) |\nabla \eta_{h, i}^\nu|_{L^2(\Omega)^N}^2
\\[3ex]
\ds
& \ds + 8 A_*^2 R_* q_1(\nu)^2 +\frac{A_* |\alpha|_{C([0, 1])}}{\delta_\alpha} \int_\Omega |\nabla v_0| \, dx.
\end{array}
\end{equation*}

We here set
\begin{equation}\label{kenApp10}
0 < \nu_* < \min \left\{ \frac{1}{32 |\alpha_0|_{C([0, 1])} A_* R_*}, ~ |\alpha|_{C([0,1])} \right\}.
\end{equation}
Therefore, applying Young's inequality again, we have
\begin{equation}\label{kenApp16-1}
\begin{array}{ll}
\multicolumn{2}{l}{
\ds \frac{A_*}{2h}(|\theta_{h, i}^\nu -v_0|_{L^2(\Omega)}^2 -|\theta_{h, i -1}^\nu -v_0|_{L^2(\Omega)}^2) -\frac{1}{4} |\nabla \eta_i|_{L^2(\Omega)^N}^2
}
\\[2ex]
& \ds +\frac{|\alpha|_{C([0, 1])}}{4\delta_\alpha} \cdot \nu |\nabla \theta_{h, i}^\nu|_{L^2(\Omega)^N}^2
 +|\alpha|_{C([0, 1])} \bigl| |\nabla \theta_{h, i}^\nu|_\nu \bigr|_{L^1(\Omega)}
\\[2ex]
\leq & \ds 10 \, R_*^2 \, q_1(\nu)^2 (1 +|v_0|_{H^1(\Omega)}^2).
\end{array}
\end{equation}

Now, taking the sum of (\ref{kenApp17}) and (\ref{kenApp16-1}) yields
\begin{equation}\label{kenApp18}
\begin{array}{ll}
\multicolumn{2}{l}
{\ds
\frac{1}{2h}(|\eta_{h, i}^\nu -w_0|_{L^2(\Omega)}^2 -|\eta_{h, i -1}^\nu -w_0|_{L^2(\Omega)}^2)
}
\\[2.5ex]
~~~~ & \ds +\frac{A_*}{2h}(|\theta_{h, i}^\nu -v_0|_{L^2(\Omega)}^2 -|\theta_{h, i -1}^\nu -v_0|_{L^2(\Omega)}^2)
\\[2.5ex]
& \ds +\frac{1}{4} |\nabla \eta_{h, i}^\nu|_{L^2(\Omega)^N}^2
 +\frac{|\alpha|_{C([0, 1])}}{4 \delta_\alpha} \cdot \nu |\nabla \theta_{h, i}^\nu|_{L^2(\Omega)^N}^2
\\[3.5ex]
& \ds +|\alpha|_{C([0, 1])} \bigl( \bigl| |\nabla \theta_{h, i}^\nu|_\nu \bigr|_{L^1(\Omega)} -\bigl| |\nabla \theta_{h, i -1}^\nu |_\nu \bigr|_{L^1(\Omega)} \bigr)
\\[2.5ex]
& \ds +\frac{\delta_\alpha}{|\alpha|_{C([0, 1])}} \int_\Omega \alpha(\eta_{h, i -1}^\nu) |\nabla \theta_{h, i -1}^\nu|_\nu \, dx
\\[3.5ex]
\multicolumn{2}{l}
{
\leq 11 q_1(\nu)^2 R_*^2 (1 +|w_0|_{H^1(\Omega)}^2 +|v_0|_{H^1(\Omega)}^2).
}
\end{array}
\end{equation}
{
Then, if we consider
$ B_* := {\delta_\alpha}/{|\alpha|_{C([0, 1])}} $ and  $ C_* := 12 R_*^2 \sup_{\nu \in (0, 1)}
 q_1(\nu)^2 $, multiplying both sides of (\ref{kenApp18}) by $ h $, we have
\begin{equation}\label{kenApp20}
\hspace{-1ex}
\begin{array}{ll}
\multicolumn{2}{l}
{
\ds \frac{1}{2}(|\eta_{h, i}^\nu -w_{0}|_{L^{2}(\Omega)}^{2} - |\eta_{h, i -1}^\nu -w_{0}|_{L^{2}(\Omega)}^{2})
}
\\[2.5ex]
& \ds + \frac{A_{\ast}}{2} (|\theta_{h, i}^\nu -v_{0}|_{L^{2}(\Omega)}^{2} - |\theta_{h, i -1}^\nu -v_{0}|_{L^{2}(\Omega)}^{2})
\\[2.5ex]
& \ds +\frac{h B_{\ast}}{4} \left( |\nabla \eta_{h, i}^\nu|_{L^{2}(\Omega)^{N}}^{2} - |\nabla \eta_{h, i -1}^\nu|_{L^{2}(\Omega)^{N}}^{2} \right)
\\[2.5ex]
& \ds + \nu \cdot \frac{h B_{\ast}}{4} \left( |\nabla \theta_{h, i}^\nu|_{L^{2}(\Omega)^{N}}^{2} - |\nabla \theta_{h, i -1}^\nu|_{L^{2}(\Omega)^{N}}^{2} \right)
\\[2.5ex]
& \ds +h |\alpha|_{C([0, 1])}(\bigl| |\nabla \theta_{h, i}^\nu|_\nu \bigr|_{L^{1}(\Omega)} - \bigl| |\nabla \theta_{h, i -1}^\nu|_\nu \bigr|_{L^{1}(\Omega)})
\\[2.5ex]
& \ds  + \frac{B_{\ast} h}{2} \mathscr{F}_{\nu}(\eta_{h, i -1}^\nu, \theta_{h, i -1}^\nu)
\\[2.5ex]
\le &  h \cdot 11 R_*^2 q_1(\nu)^2 (1 +|w_0|_{H^1(\Omega)}^2 +|v_0|_{H^1(\Omega)}^2) +h|\hat{g}|_{C([0, 1])} \mathscr{L}^N(\Omega)
\\[2.5ex]
\le & h C_*(1 +|w_0|_{H^1(\Omega)}^2 +|v_0|_{H^1(\Omega)}^2) .
\end{array}
\end{equation}
}
Thus, (\ref{lem1-inq}) is obtained  by having in mind (\hyperlink{(H3)}{H3})-(\hyperlink{(H4)}{H4}), and by taking {the sum of (\ref{kenApp20}), from the case of $i = 1$ up to that of $i = m\in\N$.}

In case that $r_1(\nu)>0$, we take instead
\begin{equation}\label{kenApp14}
A_* := \frac{2|\alpha_0|_{C([0, 1])}|\alpha|_{C([0, 1])}}{\delta_\alpha} ~ (\leq 2 \delta_\alpha R_*^{1/2}).
\end{equation}
Then, working as in the previous case, we obtain
\begin{equation*}
\begin{array}{rl}
\multicolumn{2}{l}{
\ds \frac{A_*}{2h} (|\theta_{h, i}^\nu -v_0|_{L^2(\Omega)}^2 -|\theta_{h, i -1}^\nu -v_0|_{L^2(\Omega)}^2)
}
\\[3ex]
\qquad \qquad & \ds -\left( \frac{1}{8} +4\nu |\alpha_0|_{C([0, 1])} A_* R_* \right) |\nabla \eta_{h, i}^\nu|_{L^2(\Omega)^N}^2
\\[3ex]
& \ds +2|\alpha|_{C([0, 1])} ||\nabla \theta_{h, i}^\nu|_\nu|_{L^1(\Omega)} +\frac{|\alpha|_{C([0, 1])}}{2 \delta_\alpha} \cdot \nu |\nabla \theta_{h, i}^\nu|_{L^2(\Omega)^N}^2
\\[3ex]
\ds
\leq & \ds 8 A_*^2 R_* \int_\Omega q_1(\nu)^2 |\nabla \theta_{h, i}^\nu|^{2 r_1(\nu)} \, dx
\\[3ex]
& \ds +\frac{A_* |\alpha|_{C([0, 1])}}{\delta_\alpha} \int_\Omega q_1(\nu) |\nabla v_0|^{1 +r_1(\nu)} \, dx +\frac{A_* |\alpha_0|_{C([0, 1])}}{ \delta_\alpha^2} \cdot \nu |\nabla v_0|_{L^2(\Omega)^N}^2.
\end{array}
\end{equation*}
Based on this, let us set a constant $ \nu_* \in (0, 1) $ so small to satisfy
\begin{equation}\label{kenApp10}
0 < \nu_* < \min \left\{ \frac{1}{32 |\alpha_0|_{C([0, 1])} A_* R_*}, ~ \frac{1}{2} \right\},
\end{equation}

\begin{equation}\label{kenApp11}
\left\{ ~ \parbox{6cm}{
$ \ds \frac{3}{4} \leq q_0(\nu) \leq 1 \leq q_1(\nu) \leq \frac{5}{4} $,
\\[2ex]
$ \ds 0 \leq r_0(\nu) \leq \frac{1}{4} $, $ \ds 0 \leq r_1(\nu) \leq \frac{1}{4} $,
} \right. \mbox{for any $ \nu \in (0, \nu_*) $,}
\end{equation}
and let us take a (small) constant $ \varepsilon_3 \in (0, 1) $. Then, for $ \nu \in (0, \nu_*) $, we can apply (H3)-(H4), (\hyperlink{AP2}{AP2}), (\ref{kenApp14}), (\ref{kenApp10})-(\ref{kenApp11}) and Young's inequality to deduce that, for any $0<\varepsilon_3<1$ we have
\begin{equation}\label{kenApp15}
\begin{array}{lll}
& \multicolumn{2}{l}{\ds
\frac{A_*}{2h}(|\theta_{h, i}^\nu -v_0|_{L^2(\Omega)}^2 -|\theta_{h, i -1}^\nu -v_0|_{L^2(\Omega)}^2) -\frac{1}{4} |\nabla \eta_{h, i}^\nu|_{L^2(\Omega)^N}^2}
\\[2ex]
&& \ds \quad +2 |\alpha|_{C([0, 1])} \bigl| |\nabla \theta_{h, i}^\nu)|_\nu \bigr|_{L^1(\Omega)} +\frac{|\alpha|_{C([0, 1])}}{2 \delta_\alpha} \cdot \nu |\nabla \theta_{h, i}^\nu|_{L^2(\Omega)^N}^2
\\[2ex]
\multicolumn{3}{l}{\ds
\leq \frac{25}{2} A_*^2 R_* \int_\Omega \bigl( \varepsilon_3 |\nabla \theta_{h, i}^\nu| +\varepsilon_3^{-2 r_1(\nu)/{(1 -2 r_1(\nu)})} \bigr) \, dx  +\frac{A_*}{4 \delta_\alpha} |v_0|_{H^1(\Omega)}^2}
\\[2ex]
&& \ds \qquad +\frac{5 A_* |\alpha|_{C([0, 1])}}{4 \delta_\alpha} \int_\Omega \left( \frac{1 +r_1(\nu)}{2} |\nabla v_0|^2 +\frac{1 -r_1(\nu)}{2} \right) \, dx
\\[2ex]
\multicolumn{3}{l}{\ds
\leq \varepsilon_3 \cdot \frac{50 R_*^2}{q_0(\nu)} \bigl( ||\nabla \theta_{h, i}^\nu|_\nu|_{L^1(\Omega)} +r_0(\nu) \mathscr{L}^N(\Omega) \bigr) +\frac{50 R_*^2}{\varepsilon_3} \mathscr{L}^N(\Omega)}
\\[2ex]
&& \ds \quad +\frac{25}{16} \cdot {R_*^{1/2} |\alpha|_{C([0, 1])}} \bigl( |v_0|_{H^1(\Omega)}^2 +\mathscr{L}^N(\Omega) \bigr) +\frac{R_*^{1/2}}{2} |v_0|_{H^1(\Omega)}^2
\\[2ex]
\multicolumn{3}{l}{\ds
\leq \varepsilon_3 \cdot \frac{200 R_*^2}{3} ||\nabla \theta_{h, i}^\nu|_\nu|_{L^1(\Omega)} +50 R_*^3 \left( 1 +\varepsilon_3 +\frac{1}{\varepsilon_3} \right) (1 +|v_0|_{H^1(\Omega)}^2).}
\end{array}
\end{equation}
If we choose
\begin{equation*}
\varepsilon_3 := \frac{3}{200 R_*^2} \cdot |\alpha|_{C([0, 1])} ~ (< 1),
\end{equation*}
the inequality (\ref{kenApp15}) reduces to
\begin{equation}\label{kenApp16}
\begin{array}{c}
\ds \frac{A_*}{2h}(|\theta_{h, i}^\nu -v_0|_{L^2(\Omega)}^2 -|\theta_{h, i -1}^\nu -v_0|_{L^2(\Omega)}^2) -\frac{1}{4} |\nabla \eta_i|_{L^2(\Omega)^N}^2
\\[2ex]
\ds +|\alpha|_{C([0, 1])}||\nabla \theta_{h, i}^\nu|_\nu|_{L^1(\Omega)} \leq 50 \cdot 69 R_*^{6} (1 +|v_0|_{H^1(\Omega)}^2).
\end{array}
\end{equation}

From here and letting
\begin{equation}\label{B*}
B_* := \min \left\{ \frac{1}{2} \,,\, \frac{\delta_\alpha}{|\alpha|_{C([0, 1])}} \right\} \mbox{ \ and  \ } C_* := 7 \cdot 10^3 R_*^{6},
\end{equation}
we finish the proof as in that of the previous case. \hfill $ \Box $

\section{$\Gamma$-convergence of the approximating energies}
\ \ \vspace{-3ex}

We begin this Section with the result of $ \Gamma $-convergence of the sequence of relaxed weighted-total variations (time-independent).

\begin{thm}\label{lem-gamma0}
Let us assume that
\begin{equation}\label{lem-gamma1}
\left\{
\begin{array}{ll}
\ds \beta \in W_{\rm c}(\Omega),\ \{ \beta_{\nu} \}_{\nu \in (0, 1)} \subset W_{0}(\Omega), \vspace{1mm}\\
\ds \beta_{\nu} \to \beta \ \mbox{ in } L^{2}(\Omega)\ \mbox{ and weakly in } H^{1}(\Omega)\ \mbox{ as } \nu \downarrow 0.
\end{array}
\right.
\end{equation}
Then, the sequence $\{ \Phi_{\nu}(\beta_{\nu};\ \cdot\ ) \}_{\nu \in (0, 1)}$ converges to $\Phi_{0}(\beta;\ \cdot\ )$ on $L^{2}(\Omega)$, in the sense of $\Gamma$-convergence . 
\end{thm}
\textbf{Proof.}
By (\ref{lem-gamma1}) and the strong convergence of $\{v_{\nu}\}_{\nu \in (0, 1)}$, we have
\begin{equation}\label{lem-gamma1.5}
\begin{array}{ll}
\ds \int_{\Omega} v_{\nu} \, \mbox{div}(\beta_{\nu}\varphi) \, dx &= \ds \int_{\Omega} v_{\nu} \nabla\beta_{\nu} \cdot \varphi \, dx + \int_{\Omega} v_{\nu} \beta_{\nu} \, {\rm div} \, \varphi \, dx
\\[1.5ex]
\ds &\to \ds \int_{\Omega} v \, \nabla\beta \cdot \varphi \, dx + \int_{\Omega} v \, \beta \, {\rm div} \varphi \, dx
\\[1.5ex]
\ds &= \ds \int_{\Omega} v \, {\rm div}(\beta\varphi) \, dx, \mbox{ \ as } \nu \downarrow 0,\ \mbox{ for any } \varphi \in X_{\rm c}(\Omega).
\end{array}
\end{equation}
Besides, by (\hyperlink{AP2}{AP2}), it is deduced that
\begin{equation}\label{lem-gamma2}
\begin{array}{ll}
\ds \liminf_{\nu\downarrow0} \Phi_{\nu}(\beta_{\nu};v_{\nu}) \ge \ds\liminf_{\nu\downarrow0} \int_{\Omega} \beta_{\nu} \bigl( q_{0}(\nu)|\nabla v_{\nu}|-r_{0}(\nu) \bigr) \, dx
\\[1.5ex]
\ds \ge \lim_{\nu \downarrow 0} q_{0}(\nu) \liminf_{\nu \downarrow0} \int_{\Omega} \beta_{\nu} |\nabla v_{\nu}| \, dx - \lim_{\nu\downarrow0}r_{0}(\nu)\mathscr{L}^{N}(\Omega)
\\[1.5ex]
\ds = \liminf_{\nu \downarrow 0} \Phi_{0}(\beta_{\nu} ; v_{\nu}) 
\ds \ge \lim_{\nu \downarrow 0} \int_{\Omega} v_{\nu} \, {\rm div}(\beta_{\nu}\varphi) \, dx \vspace{2mm}\\
\ds = \ds \int_{\Omega} v \, {\rm div}(\beta\varphi) \, dx,\ \mbox{ for any } \varphi \in X_{\rm c}(\Omega) \mbox{ satisfying } |\varphi| \le 1 \mbox{ a.e. in } \Omega,
\end{array}
\end{equation}
which implies the lower bound for the $\Gamma$-convergence.
\pagebreak

For the verification of the optimality, we take any $v \in D(\Phi_{0}(\beta;\ \cdot\ )) = BV(\Omega) \cap L^{2}(\Omega)$. Then, since $C^{\infty}(\overline{\Omega})$ is dense in $BV(\Omega)\cap L^{2}(\Omega)$ for the intermediate convergence, we can take a sequence $\{\omega_{n} \}_{n = 1}^\infty \subset C^{\infty}(\overline{\Omega})$, such that
\begin{equation}\label{lem1-1}
|\omega_{n} - v|_{L^{2}(\Omega)} \le 2^{-(n+1)} \mbox{ and } \bigl| |\nabla\omega_{n}|_{L^{1}(\Omega)^{N}} - |D v|(\Omega) \bigr| \le 2^{-(n+1)}, \ \mbox{ for any } n \in \N.
\end{equation}
Also, let us set 
$ p(\nu) := 1 +r_1(\nu) $ 
for any $ \nu  \in (0, 1) $,
and let us take a sequence $\{ \nu_{n} \}_{n = 1}^\infty \subset (0,1)$ such that
\begin{equation}\label{lem1-1.5}
\begin{array}{c}
\ds 0 < \nu_{n+1} < \nu_{n} < 2^{-(n+1)},\ 
\left| \int_{\Omega}\beta_{\nu}|\nabla\omega_{n}|^{p(\nu)}dx - \int_{\Omega} \beta_{\nu}|\nabla\omega_{n}| dx \right| < 2^{-(n+1)}
\\[2ex]
\ds \left| \int_{\Omega} (\beta_{\nu}q_{1}(\nu)|\nabla\omega_{n}| - \beta |\nabla\omega_{n}|) dx \right| \le 2^{-(n+1)},
\\[2ex]
\ds \mbox{and \ } \nu |\nabla\omega_{n}|_{L^{2}(\Omega)^{N}}^{2} \le 2^{-(n+1)}, \ \mbox{ for all } n \in \N \mbox{ and all } 0 < \nu \le \nu_{n}.
\end{array}
\end{equation}
Additionally, we set
\begin{equation}\label{lem1-2}
v_{\nu}:= \left\{
\begin{array}{ll}
\ds \omega_{n}, \mbox{ if } \nu_{n+1} < \nu \le \nu_{n},\ n \in \N,
\\[1ex]
\ds \omega_{1}, \mbox{ if } \nu_{1} < \nu < 1.
\end{array}
\right.
\end{equation}
Then, owing to (\ref{lem1-1})-(\ref{lem1-2}) and Remark \ref{rem-wtv} (\hyperlink{Fact1}{Fact\,1}), it is seen that
\begin{equation}\label{lem1-3}
\ds 0  \le \liminf_{\nu\downarrow0}|v_{\nu} - v|_{L^{2}(\Omega)} \le \limsup_{\nu\downarrow0}|v_{\nu}-v|_{L^{2}(\Omega)}  = \lim_{n \to \infty} |\omega_{n}-v|_{L^{2}(\Omega)} = 0,
\end{equation}
\vspace{-1ex}
\begin{equation}\label{lem1-4}
\left\{ \parbox{7cm}{
$ \ds \nabla v_{\nu} \mathscr{L}^N \to D v \ \mbox{ weakly-}\ast \mbox{ in } \M(\Omega) $,
\\[1ex]
$ \ds \int_{\Omega} \beta |\nabla v_{\nu}| dx \to \int_{\Omega} d \bigl[ \beta |D v| \bigr] $,
} \right. \mbox{as } \nu \downarrow 0.
\end{equation}

Now, taking into account (\hyperlink{AP2}{AP2}), (\ref{lem1-1})-(\ref{lem1-4}), we obtain $(\gamma2)$  as follows:
\begin{equation*}
\begin{array}{ll}
\ds 0 & \ds \le \liminf_{\nu\downarrow0}|\Phi_{\nu}(\beta_{\nu};v_{\nu}) - \Phi_{0}(\beta;v)| \le \limsup_{\nu\downarrow0}|\Phi_{\nu}(\beta_{\nu};v_{\nu}) - \Phi_{0}(\beta;v)|
\\[2ex]
\ds & \ds \le \limsup_{\nu\downarrow0} \left[ \left| \int_{\Omega} \beta_{\nu} |\nabla v_{\nu}|_{\nu} \, dx - \int_{\Omega} d[\beta |D v|] \right| + \frac{\nu}{2} | \nabla v_{\nu}|^{2}_{L^{2}(\Omega)^{N}} \right]
\\[2ex]
\ds & \ds \le \limsup_{\nu\downarrow0} \left[ \left| \int_{\Omega} \beta_{\nu} \bigl( |\nabla v_{\nu}|_{\nu}  - q_{1}(\nu) |\nabla v_{\nu}| \bigr) \, dx \right|\right.
\\[2ex]
\ds & \ds \hspace{15mm} + \left. \left| \int_{\Omega} q_{1}(\nu) \beta_{\nu} |\nabla v_{\nu}| dx - \int_{\Omega} d[\beta|D v|] \right| + \frac{\nu}{2} | \nabla v_{\nu}|^{2}_{L^{2}(\Omega)^{N}} \right]
\\[2ex]
\ds & \ds \le \lim_{\nu\downarrow0} q_{1}(\nu) \left| \int_{\Omega} \beta_{\nu}|\nabla v_{\nu}|^{p(\nu)}dx - \int_{\Omega} \beta_{\nu}|\nabla v_{\nu}|dx \right|
\\[2ex]
\ds & \ds \ \ \ + \lim_{\nu\downarrow0} \left[ |q_{1}(\nu)-q_{0}(\nu)|\int_{\Omega}\beta_{\nu}|\nabla v_{\nu}|dx + r_{0}(\nu)\L^{N}(\Omega) \right]  + \lim_{\nu\downarrow0} \frac{\nu}{2} | \nabla v_{\nu}|^{2}_{L^{2}(\Omega)^{N}}
\\[2ex]
\ds & \ds \ \ \ + \lim_{\nu\downarrow0} \left[ \left| \int_{\Omega} \bigl( q_{1}(\nu)\beta_{\nu}|\nabla v_{\nu}| - \beta |\nabla v_{\nu}| \bigr) dx \right| + \left| \int_{\Omega} \beta |\nabla v_{\nu}|dx - \int_{\Omega} d[\beta|D v|] \right| \right]
\\[2ex]
\ds & \ds = 0.
\end{array}
\end{equation*}

\hfill$\Box$
\pagebreak

For the time-dependent case, we need some auxiliary results.
{
\begin{lem}\label{lem-base}
Let $ I \subset (0, \infty) $ be a bounded open interval. Then, the following items hold:
\begin{description}
\item[\textmd{\it (\hypertarget{A-I}{\,I\,})}]If $ \beta \in \mathscr{W}_0(I; \Omega) \cap C(\overline{I}; L^2(\Omega)) $ and $ v \in L^2(I; L^2(\Omega)) $ with $ |D v({}\cdot{})|(\Omega) \in L^1(I) $, then the function $ t \in I \mapsto \Phi_0^I(\beta(t); v(t)) $ is measurable.
\item[\textmd{\it (\hypertarget{A-II}{II})}]If $ \beta \in \mathscr{W}_{\rm c}(I; \Omega) \cap C(\overline{I}; L^2(\Omega)) $ and  $ \log \beta \in  L^\infty(I \times \Omega) $, then $ \Phi_0^I(\beta;{}\cdot{}) $ is a proper, l.s.c. and convex function on $ L^2(I; L^2(\Omega)) $, and
\begin{equation*}
D(\Phi_0^I(\beta;{}\cdot{})) = \left\{ \begin{array}{l|l}
v \in L^2(I; L^2(\Omega)) & |D v({}\cdot{})|(\Omega) \in L^1(I)
\end{array} \right\}.
\end{equation*}
\end{description}
\end{lem}
}

\noindent
{
\textbf{Proof.} Let us fix $ v \in L^2(I; L^2(\Omega)) $ with $ |D v({}\cdot{})|(\Omega) \in L^1(I) $, and let us take the sequence $ \{ v_n \}_{n = 1}^\infty \subset C^\infty(\overline{I \times \Omega}) $ as in Remark \ref{rem-2} (\hyperlink{Fact6}{Fact\,4}). Then, by virtue of Remark \ref{rem-wtv} (\hyperlink{Fact1}{Fact\,1}) (cf. \cite[Lemmas 1-4]{MS}) and the assumption on $ \beta $, it is inferred that the function $ t \in I \mapsto \Phi_0(\beta(t); v(t)) $ is expressed by a limit of measurable functions $ t \in I \mapsto \Phi_0(\beta(t); v_n(t)) $, $ n \in \N $, as follows:
\begin{equation*}
\Phi_0(\beta(t); v(t)) = \lim_{n \to \infty} \Phi_0(\beta(t); v_n(t)), \mbox{ for a.e. $ t \in I $.}
\end{equation*}
Thus, item (\hyperlink{A-I}{\,I\,}) is proved. Item (\hyperlink{A-II}{II}) can be obtained as a straightforward consequence of (\hyperlink{A-I}{\,I\,}), Fatou's lemma  and the results given in (\ref{2.4}) (cf. \cite[Section 2]{MS}). \hfill $ \Box $

{
\begin{lem}\label{auxLem1}
Let $ I \subset (0, \infty) $ be a bounded open interval. Let $ \beta \in \mathscr{W}_0(I; \Omega) \cap C(\overline{I}; L^2(\Omega)) $, $ v \in L^2(I; L^2(\Omega)) $,  $ \{ \beta_n \}_{n = 1}^\infty \subset \mathscr{W}_0(I; \Omega ) \cap C(\overline{I}; L^2(\Omega)) $ and $ \{ v_n \}_{n = 1}^\infty \subset L^2(I; L^2(\Omega)) $ be such that (\ref{fact7-0}) in Remark \ref{rem-2} (\hyperlink{Fact7}{Fact\,5}) holds, and such that
\begin{equation}\label{auxLem1-00}
t \in I \mapsto \Phi_0(\beta(t); v(t)) \mbox{ and } t \in I \mapsto \Phi_0(\beta_n(t); v_n(t)), \mbox{ for $ n \in \N $, are measurable,}
\end{equation}
In addition, let us assume
\begin{equation*}
\beta \in \mathscr{W}_{\rm c}(I; \Omega), \mbox{ \ or \ } L_* := \sup_{n \in \N} \, \bigl| |D v_n({}\cdot{})|(\Omega) \bigr|_{L^1(I)} < \infty.
\end{equation*}
Then,
\begin{equation}\label{auxLem1-03}
\liminf_{n \to \infty} \int_I \int_\Omega d \bigl[ \beta_n(t)|D v_n(t)| \bigr] \, dt \geq \int_I \int_\Omega d \bigl[ \beta(t) |D v(t)| \bigr] \, dt.
\end{equation}
\end{lem}
}

\noindent
{
\textbf{Proof.}
The proof of this result is very similar to that of \cite[Lemma 6]{MS}. However, we give it for the sake of completeness. We consider only the non-trivial case: i.e.,
\begin{equation*}
\Lambda_* := \liminf_{n \to \infty} \int_I \int_\Omega d \bigl[ \beta_n(t) |D v_n(t)| \bigr] \, dt < \infty,
\end{equation*}
and we suppose that
\begin{equation*}
\int_I \int_\Omega d \bigl[ \beta_n(t) |D v_n(t)| \bigr] \, dt \to \Lambda_* \mbox{ as $ n \to \infty $,}
\end{equation*}
by taking a subsequence if necessary.
}

{
Based on this, let us first consider the case that $ \beta \in \mathscr{W}_{\rm c}(I; \Omega) $. Then, taking any $ \varphi \in X_{\rm c}(\Omega) $, we see from (\ref{fact7-0}) that
\begin{equation*}
\begin{array}{c}
\begin{array}{rcl}
\ds \int_\Omega v(t) \, {\rm div} (\beta(t) \varphi) \, dx & = & \ds \lim_{n \to \infty} \int_\Omega v_n(t) \, {\rm div} (\beta_n(t) \varphi) \, dx
\leq \liminf_{n \to \infty} \int_\Omega d [\beta_n(t) |D v_n(t)|] \, dt,
\end{array}
\\
\ \\[-1.5ex]
\mbox{for a.e. $ t \in I $, if $ |\varphi| \leq 1 $ a.e. in $ \Omega $.}
\end{array}
\end{equation*}
The above inequality implies that
\begin{equation*}
\liminf_{n \to \infty} \int_\Omega d \bigl[ \beta_n(t) |D v_n(t)| \bigr] \geq \int_\Omega d \bigl[ \beta(t) |D v(t)| \bigr], \mbox{ for a.e. $ t \in I $.}
\end{equation*}
With (\ref{auxLem1-00}) in mind, (\ref{auxLem1-03}) is deduced by using Fatou's lemma.
}
\medskip

{
Next, let us consider the case that $ L_* < \infty $. In this case,
\begin{equation}\label{auxLem1-04}
\bigl| |D v(\cdot)|(\Omega) \bigr|_{L^1(I)} \leq \liminf_{n \to \infty} \bigl| |D v_n(\cdot)|(\Omega) \bigr|_{L^1(I)} \leq L_* < \infty.
\end{equation}
Also, since $ \beta +\delta \in \mathscr{W}_{\rm c}(I; \Omega) $ and $ \{ \beta_n +\delta \}_{n = 1}^\infty \subset \mathscr{W}_{\rm c}(I; \Omega) $, for any $ \delta \in (0, 1) $, one can confirm that
\begin{equation*}
\liminf_{n \to \infty} \int_I \int_\Omega d \bigl[ (\beta_n(t) +\delta)|D v_n(t)| \bigr] \, dt \geq \int_I \int_\Omega d \bigl[ (\beta(t) +\delta)|D v(t)| \bigr] \, dt,
\end{equation*}
by applying the same argument as in the previous case of $ \beta \in \mathscr{W}_{\rm c}(I; \Omega) $.
}
\medskip

{
Hence, taking any $ \delta \in (0, 1) $, and invoking (\ref{2.4}), (\ref{auxLem1-00}) and (\ref{auxLem1-04}), it is computed that
\begin{eqnarray*}
&& \ds \liminf_{n \to \infty} \int_I \int_\Omega d \bigl[ \beta_n(t) |D v_n(t)| \bigr] \, dt
\\
& \geq & \ds \liminf_{n \to \infty} \left( \int_I \int_\Omega d \bigl[ (\beta_n(t) +\delta) |D v_n(t)| \bigr] \, dt -\delta \bigl| |D v_n(\cdot)|(\Omega) \bigr|_{L^1(I)} \right)
\\
& \geq & \ds \int_I \int_\Omega \bigl[ (\beta(t) +\delta) |D v(t)| \bigr] \, dt -\delta L_* \geq \int_I \int_\Omega \bigl[ \beta(t) |D v(t)| \bigr] \, dt -\delta L_*.
\end{eqnarray*}
Since $ \delta \in (0, 1) $ is arbitrary, the above inequality implies (\ref{auxLem1-03}). \hfill $ \Box $
}
{
\begin{rem}\label{Rem.auxCor1}
\begin{em}
As a corollary of the previous Lemma \ref{auxLem1}, we can show that if $ \beta \in W_{0}(\Omega) $, $ \{ \beta_n \}_{n = 1}^\infty \subset W_0(\Omega) $, $ v \in L^2(\Omega) $ and $ \{ v_n \}_{n = 1}^\infty \subset L^2(\Omega) $ fulfill that
\begin{equation*}
\left\{ \parbox{7.7cm}{
$ \beta_n \to \beta $ \ in $ L^2(\Omega) $ and weakly in $ H^1(\Omega) $,
\\[1ex]
$ v_n \to v $ \ in $ L^2(\Omega) $,
} \right. \mbox{as $ n \to \infty $}
\end{equation*}
and
\begin{equation*}
\beta \in W_{\rm c}(\Omega) \mbox{ \ or \ } \sup_{n \in \N} |D v_n|(\Omega) < \infty,
\end{equation*}
then,
\begin{equation*}
\liminf_{n \to \infty} \int_\Omega d \bigl[ \beta_n |D v_n| \bigr] \geq \int_\Omega d \bigl[ \beta |D v| \bigr].
\end{equation*}
\end{em}
\end{rem}
}
{
\begin{lem}\label{auxLem2}
Let $ I \subset (0, \infty) $ be a bounded open interval. Let $ \beta \in \mathscr{W}_0(I; \Omega) \cap C(\overline{I}; L^2(\Omega)) \cap L^\infty(I \times \Omega) $, $ \{ \beta_n \}_{n = 1}^\infty \subset \mathscr{W}_0(I; \Omega) \cap C(\overline{I}; L^2(\Omega)) $, $ v \in L^2(I; L^2(\Omega)) $ and $ \{ v_n \}_{n = 1}^\infty \subset $ \linebreak $ L^2(I; L^2(\Omega)) $ be such that $ |D v(\cdot)|(\Omega) $ $ \in L^1(I) $, and  (\ref{auxLem1-00}) and (\ref{fact7-0})-(\ref{fact7-1}) in Remark \ref{rem-2} (\hyperlink{Fact7}{Fact\,5}) are satisfied. Also, let $ \varrho \in C(\overline{I}; L^2(\Omega)) \cap L^\infty(I; H^1(\Omega)) \cap L^\infty(I \times \Omega) $  and $ \{ \varrho_n \}_{n = 1}^\infty \subset C(\overline{I}; L^2(\Omega)) \cap L^\infty(I; H^1(\Omega)) \cap L^\infty(I \times \Omega) $ be such that (\ref{fact7-2}) in Remark \ref{rem-2} (\hyperlink{Fact7}{Fact\,5}) holds. In addition, let us assume
\begin{equation}\label{fact7-3}
\int_I \int_\Omega d \bigl[ \beta_n(t) |D v_n(t)| \bigr] \, dt \to \int_I \int_\Omega d \bigl[ \beta(t) |D v(t)| \bigr] \, dt, \mbox{ as $ n \to \infty $.}
\end{equation}
Then,
\begin{equation}\label{auxLem2-01}
\int_I \int_\Omega d \bigl[ \varrho_n(t) |D v_n(t)| \bigr] \, dt \to \int_I \int_\Omega d \bigl[ \varrho(t) |D v(t)| \bigr] \, dt, \mbox{ as $ n \to \infty $.}
\end{equation}
\end{lem}
}

\noindent
{
\textbf{Proof.}
This lemma is proved by relying on the following elementary fact:
\begin{description}
\item[{(\hypertarget{Fact8}{Fact\,6})}]Let $ a, b \in \R $, $ \{ a_n \}_{n = 1}^\infty \subset \R $ and $ \{ b_n \}_{n = 1}^\infty \subset \R $ be such that:
\begin{equation*}
\liminf_{n \to \infty} a_n \geq a, ~ \liminf_{n \to \infty} b_n \geq b \mbox{ and } \limsup_{n \to \infty} \, (a_n +b_n) \leq a +b.
\end{equation*}
Then, $ a_n \to a $ and $ b_n \to b $ as $  n \to  \infty $.
\end{description}
}
\medskip

{
From the assumptions, it is easily checked that
\begin{equation*}
\left\{
\begin{array}{l}
\ds \{ [\varrho]^+, [\varrho_n]^+ \}_{n = 1}^\infty \subset \mathscr{W}_0(I; \Omega) \cap C(\overline{I}; L^2(\Omega)),
\\[2ex]
\left\{ \begin{array}{l|l}
\frac{M_0}{\delta_0} \beta -[\varrho]^+ , \frac{M_0}{\delta_0} \beta_n -[\varrho_n]^+ & n \in \N
\end{array} \right\} \subset  \mathscr{W}_0(I; \Omega) \cap C(\overline{I}; L^2(\Omega)).
\end{array} \right.
\end{equation*}
So, owing to (\ref{fact7-1}), we can apply Lemma \ref{auxLem1} to see that
\begin{equation}\label{auxLem2-04}
\liminf_{n \to \infty} \int_I \int_\Omega d \bigl[ [\varrho_n]^+(t) |D v_n(t)| \bigr] \, dt \geq \int_I \int_\Omega d \bigl[ [\varrho]^+(t) |D v(t)| \bigr] \, dt,
\end{equation}
together with
\begin{equation*}
\liminf_{n \to \infty} \int_I \int_\Omega d \bigl[ {\textstyle (\frac{M_0}{\delta_0} \beta_n -[\varrho_n]^+)(t)} |D v_n(t)| \bigr] \, dt \geq \int_I \int_\Omega d \bigl[ {\textstyle (\frac{M_0}{\delta_0} \beta -[\varrho]^+)(t)} |D v(t)| \bigr] \, dt.
\end{equation*}
}

{
In the meantime, from (\ref{2.4}), it follows that
\begin{equation}\label{auxLem2-02}
\begin{array}{c}
\ds \int_I \int_\Omega d \bigl[ [\varrho]^+(t)|D v(t)| \bigr] \, dt +\int_I \int_\Omega d \bigl[ {\textstyle (\frac{M_0}{\delta_0} \beta -[\varrho]^+)(t)} |D v(t)| \bigr] \, dt
\\[2ex]
= \ds \frac{M_0}{\delta_0} \int_I \int_\Omega d \bigl[ \beta(t)|D v(t)| \bigr] \, dt \leq \frac{M_0}{\delta_0} |\beta|_{L^\infty(I \times \Omega)} \bigl| |D v({}\cdot{})|(\Omega) \bigr|_{L^1(I)} < \infty.
\end{array}
\end{equation}

Therefore, taking into account (\ref{2.4}), (\ref{fact7-3}), (\ref{auxLem2-02}) and (\hyperlink{Fact3}{Fact\,3}), we can compute that
\begin{eqnarray}
&& \limsup_{n \to \infty} \left( \int_I \int_\Omega d \bigl[ [\varrho_n]^+(t)|D v_n(t)| \bigr] \, dt \right.
\nonumber
\\
&& \hspace{7ex} \left. +\int_I \int_\Omega d \bigl[ {\textstyle (\frac{M_0}{\delta_0} \beta_n -[\varrho_n]^+)(t)} |D v_n(t)| \bigr] \, dt \right)
\nonumber
\\
& = & \frac{M_0}{\delta_0} \lim_{n \to \infty} \int_I \int_\Omega d \bigl[ \beta_n(t) |D v_n(t)| \bigr] \, dt = \frac{M_0}{\delta_0}  \int_I \int_\Omega d \bigl[ \beta(t) |D v(t)| \bigr] \, dt
\label{auxLem2-03}
\\
& = & \int_I \int_\Omega d \bigl[ [\varrho]^+(t)|D v(t)| \bigr] \, dt +\int_I \int_\Omega d \bigl[ {\textstyle (\frac{M_0}{\delta_0} \beta -[\varrho]^+)(t)} |D v(t)| \bigr] \, dt.
\nonumber
\end{eqnarray}
}

{
In view of (\ref{auxLem2-04}),  (\ref{auxLem2-03}) and (\hyperlink{Fact8}{Fact\,6}), it is inferred that
\begin{equation}\label{auxLem2-05}
\lim_{n \to \infty} \int_I \int_\Omega d \bigl[ [\varrho_n]^+(t) |D v_n(t)| \bigr] \, dt = \int_I \int_\Omega d \bigl[ [\varrho]^+(t) |D v(t)| \bigr] \, dt.
\end{equation}
Similarly, one can see that
\begin{equation}\label{auxLem2-06}
\lim_{n \to \infty} \int_I \int_\Omega d \bigl[ [\varrho_n]^-(t) |D v_n(t)| \bigr] \, dt = \int_I \int_\Omega d \bigl[ [\varrho]^-(t) |D v(t)| \bigr] \, dt.
\end{equation}
The convergence (\ref{auxLem2-01}) is obtained by taking the difference between (\ref{auxLem2-05}) and (\ref{auxLem2-06}), and applying Remark \ref{rem-wtv} (\hyperlink{Fact2}{Fact\,2}). \hfill $ \Box $
}
\bigskip

{
\begin{lem}\label{auxLem10}
Let $ I \subset (0, \infty) $ be a bounded open interval. Let $ \varrho \in \mathscr{W}_{\rm c}(I; \Omega) $ and $ \{ \varrho_n \}_{n = 1}^\infty \subset \mathscr{W}_{\rm c}(I; \Omega) $ be such that
\begin{equation}\label{aL10-1}
\left\{ \parbox{12.5cm}{
	$ \log \varrho \in L^\infty(I \times \Omega) $ \ and \ $ \log \varrho_n \in L^\infty(I \times \Omega) $ for $ n \in \N $,
	\\[1ex]
	$ \varrho_n(t) \to \varrho(t) $ in $ L^2(\Omega) $ and weakly in $ H^1(\Omega) $ as $ n \to \infty $, for a.e. $ t \in I $.
} \right.
\end{equation}
Then, $ \Phi_0^I(\varrho_n;{}\cdot{}) \to \Phi_0^I(\varrho;{}\cdot{}) $ on $ L^2(I; L^2(\Omega)) $, in the sense of $ {\mit \Gamma} $-convergence, as $ n \to \infty $.
\end{lem}
}
\bigskip

\noindent
{
\textbf{Proof.}
By  assumption (\ref{aL10-1}), the lower-bound condition is a straightforward consequence of Lemmas \ref{lem-base} and \ref{auxLem1}. Also, we can verify the condition of optimality by applying Lemma \ref{auxLem2} with $ \beta = 1 $, $ \{ \beta_n \}_{n = 1}^\infty = \{ 1 \} $, and $ v \in L^2(I; L^2(\Omega)) $ and $ \{ v_n \}_{n = 1}^\infty \subset C^\infty(\overline{I \times \Omega}) $ as in Remark \ref{rem-2} (\hyperlink{Fact6}{Fact\,4}). \hfill $ \Box $
}
\bigskip

Let $ I \subset (0, \infty) $ be a fixed bounded open interval. For any $  \nu \in (0, 1) $ and any $\beta \in \mathscr{L}^{2}_{0}(I;\Omega)$, we define a functional $\Phi_{\nu}^{I}(\beta;\ \cdot\ )$ on $L^{2}(I;L^{2}(\Omega))$ as
\begin{equation} \label{phihatnu}
v \in L^{2}(I;L^{2}(\Omega)) \mapsto \Phi_{\nu}^{I}(\beta;v) := \left\{
\begin{array}{ll}
\ds \int_{I} \Phi_{\nu}(\beta(t) ; v(t) ) \, dt, \vspace{2mm}\\
\ds \hspace{0.9cm} \mbox{ if } v \in L^{2}(I;H^{1}(\Omega)), \vspace{5mm}\\
\ds \infty, \ \ \mbox{ otherwise. }
\end{array}
\right.
\end{equation}
This functional corresponds to a relaxed version of the time-dependent weighted total variation. Also, by (\hyperlink{(H3)}{H3}), it can be easily shown that the functional $\Phi_{\nu}^{I}(\beta;{}\cdot{}) $ given in (\ref{phihatnu}) is proper l.s.c. and convex in $L^{2}(I;L^{2}(\Omega))$, and $D(\Phi_{\nu}^{I}(\beta;{}\cdot{})) = L^{2}(I;H^{1}(\Omega))$ for any $ \nu \in (0, 1) $ and any $\beta \in \mathscr{L}^{2}_{0}(I;\Omega)$.

\bigskip
{Now, we show some key-properties of time-dependent weighted total variations, including the $ \Gamma $-convergence result.}

\begin{lem}[Properties kindred to the lower bound]\label{lem-gamma10}
Let $\{\nu_{n} \}_{n = 1}^\infty $ be a sequence such that $\nu_{n} \downarrow 0$ as $n \to \infty$. Also, let us assume that $\beta \in \mathscr{W}_{0}(I;\Omega) \cap C(\overline{I};L^{2}(\Omega))$, $\{\beta_{n} \}_{n = 1}^\infty \subset \mathscr{L}^{2}_{0}(I;\Omega)$, $v \in C(\overline{I};L^{2}(\Omega))$, $\{ v_{n} \}_{n = 1}^\infty \subset L^{2}(I;H^{1}(\Omega))$, and
\begin{equation}\label{lem2-assume1}
\left \{
\begin{array}{ll}
\multicolumn{2}{l}{\beta_n(t) \to \beta(t) \mbox{ in $ L^2(\Omega) $}}
\\[1ex]
& \quad \mbox{and weakly in $ H^1(\Omega) $,}
\\[1ex]
\multicolumn{2}{l}{v_n(t) \to v(t) \mbox{ in $ L^2(\Omega) $,}}
\end{array}
\right. \mbox{for a.e. $ t \in I $, as $ n \to \infty $.}
\end{equation}
In addition, let us assume that
\begin{equation}
\beta \in \mathscr{W}_{\rm c}(I;\Omega), \ \mbox{ or }\ L_{0} := \sup_{n \in \N} |\nabla v_{n}|_{L^{1}(I;L^{1}(\Omega; \R^{N}))} < \infty.
\end{equation}
Then,
\begin{equation}\label{lem2-result}
\liminf_{n \to \infty} \Phi^{I}_{\nu_{n}} (\beta_{n};v_{n}) \ge \liminf_{n \to \infty} \Phi_{0}^{I}(\beta_{n};v_{n}) \ge \Phi_{0}^{I}(\beta;v).
\end{equation}
\end{lem}
\textbf{Proof.}
We can assume that $\Lambda_{0} := \liminf_{n \to \infty} \Phi^{I}_{0}(\beta_{n};v_{n}) < \infty$, since the other case is obvious. In this case, we may suppose the existence of a subsequence, not relabeled, such that
\begin{equation}\label{lem2-1}
\Phi^{I}_{0}(\beta_{n};v_{n}) \to \Lambda_{0},\ \ \mbox{ as } n \to \infty.
\end{equation}
We first consider the case that $\beta \in \mathscr{W}_{\rm c}(I;\Omega)$. Then, {from Remark \ref{Rem.auxCor1} and (\ref{lem2-assume1})}, we can see that
\begin{equation}\label{lem2-2}
\begin{array}{ll}
\ds \liminf_{n \to \infty} \Phi_{\nu_{n}}(\beta_{n}(t);v_{n}(t))& \ds \ge \liminf_{n \to \infty} \Phi_{0}(\beta_{n}(t);v_{n}(t)) \vspace{3mm}\\
\ds & \ds \ge \Phi_{0}(\beta(t);v(t)), \ \ \mbox{ for a.e. } t \in I.
\end{array}
\end{equation}
So, by (\ref{lem2-assume1}), (\ref{lem2-1})-(\ref{lem2-2}), Remark \ref{rem-2} (\hyperlink{Fact3}{Fact\,3}), {Lemma \ref{lem-base} (\hyperlink{A-I}{\,I\,})} and Fatou's lemma, (\ref{lem2-result}) is verified as follows:
\begin{equation}\label{lem2-3}
\begin{array}{ll}
\ds \liminf_{n \to \infty} \Phi_{\nu_{n}}^{I}(\beta_{n};v_{n}) \ge \liminf_{n\to\infty} \int_{I} \Phi_{0}(\beta_{n}(t);v_{n}(t))dt = \Lambda_{0} \vspace{3mm}\\
\ds \ge \int_{I} \liminf_{n\to\infty} \Phi_{0}(\beta_{n}(t);v_{n}(t))dt \ge \int_{I}\Phi_{0}(\beta(t);v(t))dt = \Phi^{I}_{0}(\beta;v).
\end{array}
\end{equation}

Next, we consider the case of $L_{0} < \infty$. Then, it is immediately seen that
\begin{equation}\label{lem2-4}
\bigl| |D v(\cdot)|(\Omega) \bigr|_{L^{1}(I)} \le \liminf_{n \to \infty} |\nabla v_{n}|_{L^{1}(I;L^{1}(\Omega; \R^{N}))} \le L_{0} < \infty.
\end{equation}
Also, since $\beta+ \delta \in \mathscr{W}_{\rm c}(I;\Omega)$ and $\{ \beta_{n}+\delta \}_{n = 1}^\infty \subset \mathscr{W}_{\rm c}(I;\Omega)$, for any $ \delta \in (0, 1) $, we can apply the same argument as in (\ref{lem2-2})-(\ref{lem2-3}) to see that, for any $ \delta \in (0, 1) $, it holds
\begin{equation}\label{lem2-5}
\liminf_{n \to \infty} \Phi^{I}_{\nu_{n}}(\beta_{n} + \delta ;v_{n}) \ge \liminf_{n\to \infty}\Phi^{I}_{0}(\beta_{n}+\delta; v_{n}) \ge \Phi^{I}_{0}(\beta +\delta;v).
\end{equation}
On account of {(\ref{2.4})}, (\ref{lem2-assume1}), (\ref{lem2-4})-(\ref{lem2-5}) and (\hyperlink{AP2}{AP2}), it is deduced that
\begin{equation}\label{lem2-6}
\begin{array}{ll}
\ds \liminf_{n\to\infty}\Phi^{I}_{\nu_{n}}(\beta_{n};v_{n}) \ge \liminf_{n\to\infty}\Phi^{I}_{0}(\beta_{n};v_{n}) \vspace{3mm}\\
\ds = \liminf_{n \to \infty} (\Phi^{I}_{0}(\beta_{n}+\delta;v_{n}) - \delta |\nabla v_{n}|_{L^{1}(I;L^{1}(\Omega)^{N})}) \vspace{3mm}\\
\ds \ge \liminf_{n \to \infty} \Phi^{I}_{0}(\beta_{n} +\delta ; v_{n}) - \delta L_{0} \ge \Phi^{I}_{0}(\beta +\delta;v) - \delta L_{0}\vspace{3mm}\\
\ds = \int_{I}\int_{\Omega} d[(\beta(t)+\delta)|D v(t)|]dt - \delta L_{0} \vspace{3mm}\\
\ds \ge \Phi^{I}_{0}(\beta;v)-\delta (L_{0} - \bigl| |D v(\cdot)|(\Omega) \bigr|_{L^{1}(I)}),
\end{array}
\end{equation}
for any $ \delta \in (0, 1) $. Since $\delta$ is arbitrary, (\ref{lem2-6}) finishes the proof.
\hfill $\Box$

\begin{lem}[The property kindred to the optimality]\label{lem-gamma11}
Let $ I \subset (0, \infty) $ be a bounded open interval. Let $\{ \nu_{n} \}_{n = 1}^\infty \subset (0,1)$ with $\nu_{n} \downarrow 0$ as $n \to \infty$. Let $\beta \in \mathscr{W}_{0}(I;\Omega) \cap C(\overline{I};L^{2}(\Omega)) \cap L^{\infty}(I \times \Omega)$, and let $\{ \beta_{n} \}_{n = 1}^\infty \subset \mathscr{L}_{0}^{2}(I;\Omega)$ be a sequence such that $ \beta_{n}(t) \to \beta(t)$ in $L^{2}(\Omega)$ and weakly in $H^{1}(\Omega)$ as $n \to \infty$, for a.e. $t \in I$, with the additional property that there exists a constant $M_{0} > 0$ satisfying $\beta \le M_{0}$ a.e. in $I \times \Omega$, and $\beta_{n} \le M_{0}$  a.e. in $I \times \Omega$, for any $n \in \N$.
Then, for any $v \in C(\overline{I}; L^{2}(\Omega))$, satisfying $|D v(\cdot)|(\Omega) \in L^{1}(I)$, there exists a sequence $\{ v_{n} \}_{n = 1}^\infty \subset C^{\infty}(\overline{I\times\Omega})$, such that
\begin{equation}
\left \{
\begin{array}{ll}
\ds v_{n} \to v \mbox{ in } L^{2}(I ; L^{2}(\Omega)), \vspace{3mm}\\
\ds \lim_{n \to \infty} \Phi_{\nu_{n}}^{I} (\beta_{n};v_{n}) = \lim_{n \to \infty} \Phi_{0}^{I}(\beta_{n};v_{n}) = \Phi_{0}^{I}(\beta, v).
\end{array}
\right.
\end{equation}
\end{lem}
\textbf{Proof.}
By Remark \ref{rem-2} (\hyperlink{Fact6}{Fact\,6}), we can find a sequence $\{ \psi_{n} \}_{n = 1}^\infty \subset C^{\infty}(\overline{I\times\Omega})$, such that
\begin{equation}\label{lem3-1}
\psi_{n} \to v \ \mbox{ in } L^{2}(I;L^{2}(\Omega))\ \mbox{ and } \int_{I}\int_{\Omega} |\nabla\psi_{n}(t)|dxdt \to \int_{I}\int_{\Omega} d|Dv(t)|dt
\end{equation}
as $n \to \infty$. Moreover, taking a subsequence if necessary, it is possible to assume that
\begin{equation}\label{lem3-2}
\psi_{n}(t) \to v(t)\ \mbox{ in } L^{2}(\Omega), \mbox{ for a.e. } t \in I, \mbox{ as } n \to \infty.
\end{equation}
Therefore, we can apply Remark \ref{rem-2} (\hyperlink{Fact7}{Fact\,5}) with $ \varrho = 1 $ and $ \{ \varrho_n \}_{n = 1}^\infty = \{1\} $, and we obtain that
\begin{equation*}
\begin{array}{ll}
\ds \lim_{n \to \infty} \Phi_{0}^{I}(\beta_{n};\psi_{n}) & \ds = \lim_{n\to\infty}\int_{I}\int_{\Omega} \beta_{n}(t)|\nabla\psi_{n}|dxdt \vspace{3mm}\\
& \ds = \int_{I}\int_{\Omega} d[\beta(t)|D v(t)|]dt = \Phi_{0}^{I}(\beta;v).
\end{array}
\end{equation*}
In addition, we consider a sequence $\{ n_{j} \}_{j = 1}^\infty \subset \N$, such that
\begin{equation*}
\ds 2^{j} \le n_{j} < n_{j+1},\ \frac{\nu_{n}}{2}\int_{I}\int_{\Omega}|\nabla\psi_{n_{j}}|^{2}dxdt \le 2^{-(j+1)},
\end{equation*}
\begin{equation*}
\ds \left| \int_{I}\int_{\Omega} \beta_{n_{j}}(t)|\nabla\psi_{n_{j}}(t)|^{p(\nu_{n})} dxdt - \int_{I}\int_{\Omega} \beta_{n_{j}}(t)|\nabla\psi_{n_{j}}(t)|dxdt \right| < 2^{-(j+1)},
\end{equation*}
and
\begin{equation*}
\ds \left| \int_{I}\int_{\Omega} (\beta_{n_{j}}(t) q_{1}(\nu_{n})|\nabla\psi_{n_{j}}(t)| - \beta(t)|\nabla\psi_{n_{j}}(t)| ) dxdt \right| \le 2^{-(j+1)},
\end{equation*}
for any $ j \in \N $ and any $ n \geq n_j $.

Based on these, we define
\begin{equation}\label{lem3-3}
v_{n} := \left\{
\begin{array}{ll}
\ds \psi_{n_{j}} \mbox{ on } \overline{I\times\Omega}, \mbox{ if } n_{j} \le n < n_{j+1}, \vspace{3mm}\\
\ds \psi_{n_{1}} \mbox{ on } \overline{I\times\Omega}, \mbox{ if } 1 \le n < n_{1},
\end{array}
\right. \mbox{ for some } j \in \N.
\end{equation}
Then, working as in Theorem \ref{lem-gamma0} and using (\ref{lem3-1})-(\ref{lem3-3}), the desired convergence is obtained.
\hfill $\Box$

{
\begin{thm}\label{lem-gamma}
Let $ I \subset (0, \infty) $ be a bounded open interval. Let $\{ \nu_{n} \}_{n = 1}^\infty \subset (0,1)$ with $\nu_{n} \downarrow 0$ as $n \to \infty$. Let $ \beta \in \mathscr{W}_{\rm c}(I; \Omega) $ and $ \{ \beta_n \}_{n = 1}^\infty \subset \mathscr{W}_{\rm c}(I; \Omega) $ be as in (\ref{aL10-1}).
Then, $ \Phi_0^I(\beta_n;{}\cdot{}) \to \Phi_0^I(\beta;{}\cdot{}) $ on $ L^2(I; L^2(\Omega)) $, in the sense of $ {\mit \Gamma} $-convergence, as $ n \to \infty $.
\end{thm}
}

\noindent
{
\textbf{Proof. }
By virtue of the assumption (\ref{aL10-1}), this result is obtained as a straightforward consequence of Lemmas \ref{lem-gamma10}-\ref{lem-gamma11}. \hfill $ \Box $}

{
\begin{lem}\label{Lem.A_last}
Let $ I \subset (0, \infty) $ be a bounded open interval. Let $\{ \nu_{n} \}_{n = 1}^\infty \subset (0,1)$ with $\nu_{n} \downarrow 0$ as $n \to \infty$. Let us assume that $ \beta \in \mathscr{W}_0(I; \Omega) \cap C(\overline{I}; L^2(\Omega)) \cap L^\infty(I \times \Omega) $, $ \{ \beta_n \}_{n = 1}^\infty \subset \mathscr{L}_0^2(I; \Omega) $, $ v \in L^2(I; L^2(\Omega)) $ with $ |Dv({}\cdot{})|(\Omega) \in L^1(I) $, $ \{ v_n \}_{n = 1}^\infty \subset L^2(I; H^1(\Omega)) $, $ \varrho \in C(\overline{I}; L^2(\Omega)) \cap L^\infty(I; H^1(\Omega)) \cap L^\infty(I \times \Omega) $ and $ \{ \varrho_n \}_{n = 1}^\infty \subset L^2(I; L^2(\Omega)) $, and these fulfill (\ref{fact7-0})-(\ref{fact7-2}) in Remark \ref{rem-2} (\hyperlink{Fact7}{Fact\,5}). In addition, let us assume that
\begin{equation*}
\ds \int_I \int_\Omega \beta_n(t) |\nabla v_n(t)|_{\nu_n} \, dx dt \to \int_I \int_\Omega d \bigl[ \beta(t) |Dv(t)| \bigr] \, dt, \mbox{as $ n \to \infty $,}
\end{equation*}
then
\begin{equation*}
\int_I \int_\Omega \varrho_n(t) |\nabla v_n(t)|_{\nu_n} \, dx dt \to \int_I \int_\Omega d \bigl[ \varrho(t) |D v(t)| \bigr] \, dt \mbox{ \ as $ n \to \infty $.}
\end{equation*}
\end{lem}
}

\noindent
{
\textbf{Proof. }
The proof of this lemma is a slight modification of that of Lemma \ref{auxLem2} (or \cite[Lemma 7]{MS}) and we omit it. \hfill $ \Box $
}

\section{Proof of Main Theorem 1}
\ \ \vspace{-3ex}

This Section is devoted to the proof of Main Theorem 1. Let $ \nu_*, h_* \in (0, 1) $ be the constants as in Theorem \ref{th1} and Lemma \ref{lem1}. Let $ \{ [\tilde{\eta}_0^\nu, \tilde{\theta}_0^\nu] \}_{\nu \in (0, \nu_*)} $ be  such that
\begin{equation}\label{kenV-01}
\left\{ \begin{array}{l}
[\tilde{\eta}_0^\nu, \tilde{\theta}_0^\nu] \in D_*(\theta_0), \mbox{ \ for all $ \nu \in (0, \nu_*) $,}
\\[1ex]
[\tilde{\eta}_0^\nu, \tilde{\theta}_0^\nu] \to [\eta_0, \theta_0] \mbox{ in $ L^2(\Omega)^2 $, \ as $ \nu \downarrow 0 $.}
\end{array} \right.
\end{equation}

\medskip

{For any $ \nu \in (0, \nu_*) $ and $ h  \in (0, h_*) $,} we denote by $ \{ [\tilde{\eta}_{h, i}^\nu, \tilde{\theta}_{h, i}^\nu] \}_{i = 1}^\infty $ the unique solution to $ \mbox{(AP$_h^\nu $)} $ in the case when $ [{\eta}_0^\nu, \theta_0^\nu] = [\tilde{\eta}_0^\nu, \tilde{\theta}_0^\nu] $, provided by Theorem \ref{th1} and we define three different kinds of time-interpolations $ [\overline{\eta}_h^\nu, \overline{\theta}_h^\nu] \in L_{\rm loc}^\infty([0, \infty); H^1(\Omega))^2 $, $ [\underline{\eta}_h^\nu, \underline{\theta}_h^\nu] \in L_{\rm loc}^\infty([0, \infty); H^1(\Omega))^2 $ and $ [\widehat{\eta}_h^\nu, \widehat{\theta}_h^\nu] \in W_{\rm loc}^{1, \infty}([0, \infty); H^1(\Omega))^2 $, by letting
\begin{equation}\label{kenV-20}
\left \{
\begin{array}{ll}
\ds [\overline{\eta}_{h}^{\nu}(t), \overline{\theta}_{h}^{\nu}(t)] := [\tilde{\eta}_{h,i}^{\nu}, \tilde{\theta}_{h,i}^{\nu}], & \mbox{if $ t \in ((i -1)h, ih] \cap [0, \infty) $ with some $ i \in \mathbb{Z} $,}
\\[2ex]
[\underline{\eta}_{h}^{\nu}(t), \underline{\theta}_{h}^{\nu}(t)] := [\tilde{\eta}_{h,i-1}^{\nu}, \tilde{\theta}_{h,i-1}^{\nu}], & \mbox{if $ t \in [(i -1)h, ih) $ with some $ i \in \N $,}
\\[1ex]
\multicolumn{2}{l}{
\ds [\widehat{\eta}_{h}^{\nu,n}(t), \widehat{\theta}_{h}^{\nu,n}(t)] := \frac{ih -t}{h}[\tilde{\eta}_{h,i-1}^{\nu}, \tilde{\theta}_{h,i-1}^{\nu}] + \frac{t-(i -1)h}{h}[\tilde{\eta}_{h,i}^{\nu}, \tilde{\theta}_{h,i}^{\nu}],
}
\\[1ex]
& \mbox{if $ t \in [(i -1)h, ih) $ with some $ i \in \N $,}
\end{array}
\right.
\end{equation}
for all $ t \geq 0 $. Then, from (\ref{kenV-01}), we immediately see that
\begin{equation}\label{app-rist}
\begin{array}{c}
\left\{
\begin{array}{l}
\ds 0 \leq \overline{\eta}_h^\nu(t) \leq 1, ~~ 0 \leq \underline{\eta}_h^\nu(t) \leq 1 \mbox{ \ and \ } 0 \leq \widehat{\eta}_h^\nu(t) \leq 1,
\\[0.5ex]
\ds \max \left\{ |\overline{\theta}_h^\nu(t)|, ~ |\underline{\theta}_h^\nu(t)|, ~ |\widehat{\theta}_h^\nu(t)| \right\} \leq |\tilde{\theta}_0^\nu|_{L^\infty(\Omega)} \leq |\theta_0|_{L^\infty(\Omega)},
\end{array}
\right. \mbox{a.e in $ \Omega $,}
\ \\
\\[-2.5ex]
\mbox{for all $ t \geq 0 $, $ \nu \in (0, \nu_*) $ and $ h \in (0, h_*) $.}
\end{array}
\end{equation}
Additionally, from the inequalities (\ref{ene-inq}) and (\ref{lem2-inq}) in Theorem \ref{th1} and (\ref{lem1-inq}) in Lemma \ref{lem1}, it is inferred that
\begin{equation}\label{kenV-100}
\begin{array}{c}
\ds \frac{1}{2} \int_s^t |(\widehat{\eta}_h^\nu)_t(\tau)|_{L^2(\Omega)}^2 \, d \tau +\int_s^t \left| {\textstyle \sqrt{\alpha_0(\overline{\eta}_h^\nu(\tau))}} (\widehat{\theta}_h^\nu)_t(\tau) \right|_{L^2(\Omega)}^2 \, d \tau
\\[1.5ex]
+\mathscr{F}_{\nu_n}(\overline{\eta}_h^\nu(t), \overline{\theta}_h^\nu(t))  \leq \mathscr{F}_{\nu}(\underline{\eta}_h^\nu(s), \underline{\theta}_h^\nu(s)), 
\end{array}
\end{equation}
\begin{equation}\label{kenV-04}
\begin{array}{rl}
\ds \frac{1}{2} \int_0^t \tau |(\widehat{\eta}_h^\nu)_t(\tau)|_{L^2(\Omega)}^2 \, d \tau & \ds + ~ \int_0^t \tau |{\textstyle \sqrt{\alpha_0(\overline{\eta}_h^\nu(\tau))}}(\widehat{\theta}_h^\nu)_t(\tau)|_{L^2(\Omega)} \, d \tau
\\[1.5ex]
& \ds + ~ t \mathscr{F}_\nu(\overline{\eta}_h^\nu(t), \overline{\theta}_h^\nu(t)) \leq \int_0^t \mathscr{F}_\nu(\underline{\eta}_h^\nu(\tau), \underline{\theta}_h^\nu(\tau)) \, d \tau,
\end{array}
\vspace{-1ex}
\end{equation}
and
\vspace{-1ex}
\begin{equation}\label{kenV-03}
\begin{array}{c}
\ds \frac{1}{2} \left( |\overline{\eta}_h^\nu(t) -w_0|_{L^2(\Omega)}^2 +A_* |\overline{\theta}_h^\nu(t) -v_0|_{L^2(\Omega)}^2 \right) +{\frac{B_*}{2}} \int_0^t \mathscr{F}_\nu(\underline{\eta}_h^\nu(\tau), \underline{\theta}_h^\nu(t)) \, d \tau
\\[1.5ex]
\leq \ds \frac{1}{2} \left( |\tilde{\eta}_0^\nu -w_0|_{L^2(\Omega)}^2 +A_* |\tilde{\theta}_0^\nu -v_0|_{L^2(\Omega)}^2 \right) +\frac{h}{B_*} \mathscr{F}_\nu(\tilde{\eta}_0^\nu, \tilde{\theta}_0^\nu)
\\[1.5ex]
\ds +2 t C_*(1 +|w_0|_{H^1(\Omega)}^2 +|v_0|_{H^1(\Omega)}^2),
\end{array}
\vspace{-1ex}
\end{equation}
for all $ 0 \leq s \leq t < \infty $, $ \nu \in (0, \nu_*) $ and $ h \in (0, h_*) $, respectively.

\medskip
Now, by a diagonal type argument, we can obtain sequences $ \{ \nu_n \}_{n = 1}^\infty \subset (0, \nu_*) $ and $ \{ h_n \}_{n = 1}^\infty \subset (0, h_*) $, such that
\begin{equation}\label{kenV-02}
\left\{ ~ \parbox{9cm}{
$ \ds 0 < \nu_{n +1} < \nu_n < \nu_* 2^{-n} $, ~ $ \ds 0 < h_{n +1} < h_n < h_* 2^{-n} $,
\\[1ex]
$ \ds 0 \leq h_n \mathscr{F}_{\nu_n}(\tilde{\eta}_0^{\nu_n}, \tilde{\theta}_0^{\nu_n}) < 2^{-n} $,
} \right. \mbox{for all $ n \in \N $.}
\end{equation}
Due to (\ref{app-rist})-(\ref{kenV-02}), the sequences
\begin{equation}\label{kenV-29}
\left\{ \parbox{14cm}{
$ \{ [\overline{\eta}_n, \overline{\theta}_n] \}_{n = 1}^\infty := \{ [\overline{\eta}_{h_n}^{\nu_n}, \overline{\theta}_{h_n}^{\nu_n}] \}_{n = 1}^\infty $, \ $ \{ [\underline{\eta}_n, \underline{\theta}_n] \}_{n = 1}^\infty := \{ [\underline{\eta}_{h_n}^{\nu_n}, \underline{\theta}_{h_n}^{\nu_n}] \}_{n = 1}^\infty $,
\\[1ex]
$ \{ [\widehat{\eta}_n, \widehat{\theta}_n] \}_{n = 1}^\infty := \{ [\widehat{\eta}_{h_n}^{\nu_n}, \widehat{\theta}_{h_n}^{\nu_n}] \}_{n = 1}^\infty $, \ and \ $ \{ [\eta_{0, n}, \theta_{0, n}] \}_{n = 1}^\infty := \{ [\tilde{\eta}_0^{\nu_n}, \theta_0^{\nu_n}] \}_{n = 1}^\infty $,
} \right.
\vspace{0.5ex}
\end{equation}
satisfy the following properties:
\vspace{0.5ex}
\begin{description}
\item[\textmd{(\hypertarget{1-a}{$\sharp$1-a})}]$ \left\{
[\overline{\eta}_n(t), \overline{\theta}_n(t)], ~ [\underline{\eta}_n(t), \underline{\theta}_n(t)], ~ [\widehat{\eta}_n(t), \widehat{\theta}_n(t)] \right\}_{n = 1}^\infty \subset D_*(\theta_0) $, \ for all $ t \geq 0 $;
\vspace{0.5ex}
\item[\textmd{($\hypertarget{1-b}{\sharp$1-b})}]$ \{ [\overline{\eta}_n, \overline{\theta}_n] \}_{n = 1}^\infty $ and $ \{ [\underline{\eta}_n, \underline{\theta}_n] \}_{n = 1}^\infty $ are bounded in $ L_{\rm loc}^\infty((0, \infty); L^2(\Omega))^2 $, and $ \{ [\widehat{\eta}_n, \widehat{\theta}_n] \}_{n = 1}^\infty $ is bounded in $ C_{\rm loc}((0, \infty); L^2(\Omega))^2 $;
\vspace{0.5ex}
\item[\textmd{(\hypertarget{1-c}{$\sharp$1-c})}]$ \{ \mathscr{F}_{\nu_n}(\overline{\eta}_n({}\cdot{}), \overline{\theta}_n({}\cdot{})) \}_{n = 1}^\infty $ and $ \{ \mathscr{F}_{\nu_n}(\underline{\eta}_n({}\cdot{}), \underline{\theta}_n({}\cdot{})) \}_{n = 1}^\infty $ are sequences of nonincreasing functions on $ (0, \infty) $, which are bounded in $ L_{\rm loc}^1([0, \infty)) $ 
 and $ BV_{\rm loc}((0, \infty)) $;
\item[\textmd{(\hypertarget{1-d}{$\sharp$1-d})}]$ h_n \mathscr{F}_{\nu_n}(\eta_{0, n}, \theta_{0, n}) \to 0 $ \ as $ n \to \infty $.
\end{description}
Therefore, by {the compactness theories as in \cite{AFP,Simon}}, we can find  $ [\eta, \theta] \in C_{\rm loc}((0, \infty); L^2(\Omega))^2 $ and a function $ {\mathscr{J}}_* \in BV_{\rm loc}((0, \infty)) $, together with subsequences (not relabeled) of
\linebreak
$ \{ [\overline{\eta}_n, \overline{\theta}_n] \}_{n = 1}^\infty $, $ \{ [\underline{\eta}_n, \underline{\theta}_n] \}_{n = 1}^\infty $ and $ \{ [\widehat{\eta}_n, \widehat{\theta}_n] \}_{n = 1}^\infty $ , such that
\begin{equation}\label{kenV-05}
\left\{ \hspace{-2ex} \parbox{10.5cm}{
\vspace{-2ex}
\begin{itemize}
\item $ \eta \in W_{\rm loc}^{1, 2}((0, \infty); L^2(\Omega)) \cap L_{\rm loc}^\infty((0, \infty); H^1(\Omega)) $, and $ 0 \leq \eta \leq 1 $ a.e. in $ Q $,
\item $ \theta \in W_{\rm loc}^{1, 2}((0, \infty); L^2(\Omega)) $, $ |D \theta(\cdot)|(\Omega) \in L_{\rm loc}^\infty((0, \infty)) $, and $ |\theta| \leq |\theta_0|_{L^\infty(\Omega)} $ a.e. in $ Q $;
\vspace{-2ex}
\end{itemize}
} \right.
\end{equation}
\begin{equation}\label{kenV-06}
\left\{ \hspace{-2ex} \parbox{11cm}{
\vspace{-2ex}
\begin{itemize}
\item $ \overline{\eta}_n \to \eta $ and $ \underline{\eta}_n \to \eta $ in $ L_{\rm loc}^\infty((0, \infty); L^2(\Omega)) $, weakly-$ * $ in $ L_{\rm loc}^\infty((0, \infty); H^1(\Omega)) $, and weakly-$*$ in $ L^\infty(Q) $,
\item $ \widehat{\eta}_n \to \eta $ in $ C_{\rm loc}((0, \infty); L^2(\Omega)) $, weakly-$ * $ in $ L_{\rm loc}^\infty((0, \infty); H^1(\Omega)) $, and weakly-$*$ in $ L^\infty(Q) $,
\item $ \overline{\eta}_n(t) \to \eta(t) $, $ \underline{\eta}_n(t) \to \eta(t) $ and $ \widehat{\eta}_n(t) \to \eta(t) $ in $ L^2(\Omega) $, weakly in $ H^1(\Omega) $, for any $ t > 0 $,
\vspace{-2ex}
\end{itemize}
} \right. \mbox{ \ \ \ \ as $ n \to \infty $;}
\end{equation}
\begin{equation}\label{kenV-07}
\left\{ \hspace{-2ex} \parbox{11cm}{
\vspace{-2ex}
\begin{itemize}
\item $ \overline{\theta}_n \to \theta $ and $ \underline{\theta}_n \to \theta $ in $ L_{\rm loc}^\infty((0, \infty); L^2(\Omega)) $, and weakly-$*$ in $ L^\infty(Q) $,
\item $ \widehat{\theta}_n \to \theta $ in $ C_{\rm loc}((0, \infty); L^2(\Omega)) $,  and weakly-$*$ in $ L^\infty(Q) $,
\item $ \overline{\theta}_n(t) \to \theta(t)  $, $ \underline{\theta}_n(t) \to \theta(t) $ and $ \widehat{\theta}_n(t) \to \theta(t) $ in $ L^2(\Omega) $, and weakly-$ * $ in $ BV(\Omega) $, for any $ t > 0 $,
\vspace{-2ex}
\end{itemize}
} \right. \mbox{ \ \ \ \ as $ n \to \infty $;}
\end{equation}
and
\begin{equation}\label{kenV-08}
\begin{array}{c}
\ds \mathscr{F}_{\nu_n}(\underline{\eta}_n, \underline{\theta}_n) \to {\mathscr{J}}_*
\mbox{ \ weakly-$*$ in $ BV_{\rm loc}((0, \infty)) $, weakly-$*$ in $ L_{\rm loc}^\infty((0, \infty)) $,}
\\[1ex]
\mbox{and a.e. in $ (0, \infty) $, \ as $ n \to \infty $.}
\end{array}
\end{equation}

Next, we will prove that the limit $ [\eta, \theta] $ satisfies the variational inequalities (\ref{vari01}) and (\ref{vari02}). Let $ I $ be any bounded open interval such that  $ I \subset\subset (0, \infty) $. Then, by (\ref{kenApp01}) and (\ref{kenApp02}), the sequences in (\ref{kenV-05})-(\ref{kenV-07}) fulfill the following two variational formulas:
\begin{equation}\label{kenV-09}
\begin{array}{c}
\ds \int_I \bigl( (\widehat{\eta}_n)_t(t) +g(\overline{\eta}_n(t)), \omega(t) \bigr)_{L^2(\Omega)} \, dt +\int_I \bigl( \nabla \overline{\eta}_n(t), \nabla \omega(t) \bigr)_{L^2(\Omega)^N} \, dt
\\[1ex]
\ds +\int_I \int_\Omega \omega(t) \alpha'(\overline{\eta}_n(t)) {|\nabla \underline{\theta}_n(t)|_{\nu_n}} \, dx dt = 0,
\\[2ex]
\mbox{for any $ \omega \in L^2(I; H^1(\Omega)) \cap L^\infty(I \times \Omega) $ and any $ n \in \N $,}
\end{array}
\end{equation}
and
\begin{equation}\label{kenV-10}
\begin{array}{l}
\ds \int_I \bigl( \alpha_0(\overline{\eta}_n(t)) (\widehat{\theta}_n)_t(t), \overline{\theta}_n(t) -v(t)  \bigr)_{L^2(\Omega)} \, dt +\Phi_{\nu_n}^I(\alpha(\overline{\eta}_n); \overline{\theta}_n)
\\[1ex]
\qquad \ds \leq \Phi_{\nu_n}^I(\alpha(\overline{\eta}_n); v),
\mbox{ \ for any $ v \in L^2(I; H^1(\Omega)) $ and any $ n \in \N $.}
\end{array}
\end{equation}

{
Let us take any $ v \in C(\overline{I}; L^2(\Omega)) $ with $ |Dv({}\cdot{})|(\Omega) \in L^1(I) $. Then, on account of (\ref{2.4}), (\ref{2.4-1}) and (\ref{kenV-05}), it follows that $ v \in D(\Phi_0^I(\alpha(\eta));{}\cdot{})) $. So, by Lemma \ref{lem-gamma}}, we can find a sequence $ \{ v_n \}_{n = 1}^\infty \subset L^2(I; H^1(\Omega)) $, such that
\begin{equation*}
v_n \to v \mbox{ in $ L^2(I; L^2(\Omega)) $ and $ \Phi_{\nu_n}^I(\alpha(\overline{\eta}_n); v_n) \to \Phi_0^I(\alpha(\eta); v) $, \ as $ n \to \infty $.}
\end{equation*}
Then, by (\ref{kenV-05})-(\ref{kenV-07}) and Lemma \ref{lem-gamma}, letting $ n \to \infty $ in (\ref{kenV-10}), we obtain
{
\begin{equation}\label{kenV-11}
\begin{array}{ll}
\multicolumn{2}{l}{\ds \int_{I} (\alpha_{0}(\eta(t)) \theta_{t}(t), \theta(t) -v(t))_{L^{2}(\Omega)} dt + \Phi_{0}^{I}(\alpha(\eta);\theta)}
\\
~~~~ & \ds \le \lim_{n \to \infty} \int_{I} (\alpha_{0}(\overline{\eta}_{n})(\widehat{\theta}_{n})_{t}(t), \overline{\theta}_{n}(t)-v_n(t))_{L^{2}(\Omega)} dt + \liminf_{n\to\infty} \Phi_{{\nu}_{n}}^{I}(\alpha(\overline{\eta}_{n}); \overline{\theta}_{n})
\\
& \ds \le \int_{I} (\alpha_{0}(\eta) \theta_{t}(t), \theta(t)-v(t))_{L^{2}(\Omega)} dt + \limsup_{n\to\infty} \Phi_{{\nu}_{n}}^{I}(\alpha(\overline{\eta}_{n}); \overline{\theta}_{n})
\\
& \ds \le \lim_{n \to \infty} \Phi_{{\nu}_{n}}^{I}(\alpha(\overline{\eta}_{n});v_n)= \Phi_{0}^{I}(\alpha(\eta);v),
\\
\multicolumn{2}{l}{\mbox{ for any $ v \in C(\overline{I}; L^2(\Omega)) $ with $ |Dv({}\cdot{})|(\Omega) \in L^1(I) $.}}
\end{array}
\end{equation}
}
Since the choice of the bounded open interval $ I \subset \hspace{-0.25ex} \subset (0, \infty) $ is arbitrary, we can see from (\ref{kenV-11}) that the pair $ [\eta, \theta] $ fulfills the variational inequality (\ref{vari02}).
\medskip

{
Next, letting $ v = \theta $ in (\ref{kenV-11}), it is observed that
\begin{equation}\label{kenV-17}
\lim_{n \to \infty} \Phi_{\nu_n}^I(\alpha(\overline{\eta}_n); \overline{\theta}_n) = \Phi_0^I(\alpha(\eta); \theta),
\end{equation}
and accordingly,
\begin{equation}\label{kenV-14}
\begin{array}{rcl}
0 & \leq & \ds \liminf_{n \to \infty} \frac{\nu_n}{2} \int_I \int_\Omega |\nabla \overline{\theta}_n(t)|^2 \, dx dt \leq \limsup_{n \to \infty} \frac{\nu_n}{2} \int_I \int_\Omega |\nabla \overline{\theta}_n|^2 \, dx dt
\\[2ex]
& \leq & \ds \lim_{n \to \infty} \Phi_{\nu_n}^I(\alpha(\overline{\eta}_n); \overline{\theta}_n) -\liminf_{n \to \infty} \Phi_0^I(\alpha(\overline{\eta}_n); \overline{\theta}_n) \leq 0.
\end{array}
\end{equation}
}
In particular, (\ref{kenV-17}) and (\ref{kenV-14}) imply that
\begin{equation}\label{kenV-18}
\int_I \int_\Omega \alpha(\overline{\eta}_n(t)) {|\nabla \overline{\theta}_n(t)|_{\nu_n}} \, dx dt \to \int_I \int_\Omega d \bigl[ \alpha(\eta(t)) |D \theta(t)| \bigr] \, dt \mbox{ \ as $ n \to \infty $.}
\end{equation}

Having in mind (\ref{kenV-05})-(\ref{kenV-07}) and (\ref{kenV-18}), we apply Lemma \ref{Lem.A_last} with
$ \beta = \alpha(\eta) $, $ \{ \beta_n \}_{n = 1}^\infty = \{ \alpha(\overline{\eta}_n) \}_{n = 1}^\infty $, $ v = \theta $, $ \{ v_n \}_{n = 1}^\infty = \{ \overline{\theta}_n \}_{n = 1}^\infty $, $ \varrho = 1 $ and $ \{ \varrho_n \}_{n = 1}^\infty = \{ 1 \} $. Then,
{
\begin{equation}\label{kenV-34}
\int_I \int_\Omega |\nabla \overline{\theta}_n(t)|_{\nu_n} \, dx dt \to \int_I \int_\Omega |D \theta(t)| \, dt \mbox{ \ as $ n \to \infty $.}
\end{equation}
}
On the other hand, from (H4), \eqref{ene-inq}, (\ref{kenV-20}) and (\ref{kenV-02}), we infer that
{
\begin{equation}\label{kenV-35}
\ds \left| \int_{I}\int_{\Omega} \bigl( |\nabla \underline{\theta}_{n}|_{\nu_n} -|\nabla \overline{\theta}_{n}|_{\nu_n} \bigr) \, dxdt \right| \le \frac{2}{\delta_\alpha} h_n \mathscr{F}_{\nu_n}(\eta_{0, n}, \theta_{0, n}) \to 0,
\mbox{ as $ n \to \infty $.}
\end{equation}
}
Hence, taking any $w  \in H^{1}(\Omega) \cap L^{\infty}(\Omega) $, and applying {Lemma \ref{Lem.A_last}} with $ \beta = 1 $, $ \{ \beta_{n} \}_{n = 1}^\infty = \{ 1 \} $, $ v = \theta  $, $ \{ v_n \}_{n = 1}^\infty = \{ \underline{\theta}_n \}_{n = 1}^\infty $, $ \varrho = w \alpha'(\eta) $ and $ \{ \varrho_{n} \}_{n = 1}^\infty = \{ w \alpha'(\overline{\eta}_{n}) \}_{n = 1}^\infty $, we obtain
{
\begin{equation} \label{1steq-meas}
\begin{array}{c}
\ds \int_{I} \int_{\Omega} w \alpha'(\overline{\eta}_{n}(t))|\nabla \underline{\theta}_{n}(t)|_{\nu_n} dxdt \to \int_{I} \int_{\Omega} d \bigl[ w \alpha'(\eta(t)) |D \theta(t)| \bigr] \, dt
\\[2ex]
\mbox{as $ n \to \infty $, \ for any $ w \in H^1(\Omega) \cap L^\infty(\Omega) $.}
\end{array}
\end{equation}
}

On account of (\ref{kenV-05})-(\ref{kenV-07}) and (\ref{1steq-meas}), letting $n \to \infty$ in (\ref{kenV-09}) yields
\begin{equation}\label{kenV-22}
\begin{array}{c}
\ds \int_{I} (\eta_{t}(t) + g(\eta(t)), w)_{L^{2}(\Omega)} dt + \int_{I}\bigl( \nabla \eta(t), \nabla w \bigr)_{L^2(\Omega)^N} \, dxdt
\\[2ex]
\ds + \int_{I}\int_{\Omega} d \bigl[ w \alpha'(\eta(t))|D\theta(t)| \bigr] \, dt = 0, \mbox{ \ for any $w \in H^{1}(\Omega) \cap L^{\infty}(\Omega)$. }
\end{array}
\end{equation}
Since $I$ is arbitrary, $ [\eta, \theta] $ satisfies the variational inequality (\ref{vari01}).
\bigskip

In order to finish the proof, it remains to prove the following three items:

\begin{description}
\item[\textmd{(\hypertarget{2-a}{$\sharp$2-a})}]$ \eta \in L_{\rm loc}^2([0, \infty); H^1(\Omega)) $ and $ |D \theta({}\cdot{})|(\Omega) \in L_{\rm loc}^1([0, \infty)) $;
\vspace{0.5ex}
\item[\textmd{(\hypertarget{2-b}{$\sharp$2-b})}]$ [\eta, \theta] \in C([0, \infty); L^2(\Omega))^2 $ and $ [\eta(0), \theta(0)] = [\eta_0, \theta_0] $ in $ L^2(\Omega)^2 $;
\item[\textmd{(\hypertarget{2-c}{$\sharp$2-c})}]$ \mathscr{J}_* $ obtained in (\ref{kenV-08}) is nonincreasing, and $ \mathscr{J}_* = \mathscr{F}(\eta, \theta) $ a.e on $ (0, \infty) $.
\end{description}

We fix $ t \in (0, \infty) $, $ n \in \N $ and $ \ell \in \N $, and consider the inequality (\ref{kenV-03}) with $ [h, \nu] = [h_n, \nu_n] $ and $ [w_0, v_0] = [\eta_{0, \ell}, \theta_{0, \ell}] $. Then, with (\hyperlink{(H4)}{H4}), (\ref{relaxFreeEnrg}), (\ref{kenV-01}), (\ref{kenV-02})-(\ref{kenV-07}) and Theorem \ref{lem-gamma} in mind, letting $ n \to \infty $ yields that
\begin{eqnarray}
\frac{1}{2} \bigl( |\eta(t) -\eta_{0, \ell}|_{L^2(\Omega)}^2 +A_* |\theta(t) -\theta_{0, \ell}|_{L^2(\Omega)}^2 \bigr)
+\frac{1}{2} |\nabla \eta|_{L^2(0, t; L^2(\Omega)^N)}^2 +\delta_\alpha \bigl| |D \theta({}\cdot{})|(\Omega) \bigr|_{L^1(0, t)}
\nonumber
\\ \nonumber
 \leq \frac{1}{2} \lim_{n \to \infty} \bigl( |\overline{\eta}_n(t) -\eta_{0, \ell}|_{L^2(\Omega)}^2 +A_* |\overline{\theta}_n(t) -\theta_{0, \ell}|_{L^2(\Omega)}^2 \bigr)
+\frac{B_*}{2} \liminf_{n \to \infty} \int_0^t \mathscr{F}_{\nu_n}(\underline{\eta}_n(\tau), \underline{\theta}_n(\tau)) \, d \tau
\label{kenV-30}
\\
 \leq \frac{1}{2} \bigl( |\eta_{0, \ell} -\eta_{0}|_{L^2(\Omega)}^2 +A_* |\theta_{0, \ell} -\theta_{0}|_{L^2(\Omega)}^2 \bigr)
+2t C_* \bigl( 1 +|\eta_{0, \ell}|_{H^1(\Omega)}^2 +|\theta_{0, \ell}|_{H^1(\Omega)}^2 \bigr),\quad
\end{eqnarray}
which, together with (\ref{kenV-05}), implies (\hyperlink{2-a}{$\sharp$2-a}).
\medskip

In the meantime, since
\begin{equation*}
\begin{array}{c}
\ds |\eta(t) -\eta_0|_{L^2(\Omega)}^2 +A_* |\theta(t) -\theta_0|_{L^2(\Omega)}^2 \leq 2 \bigl( |\eta(t) -\eta_{0, \ell}|_{L^2(\Omega)}^2 +A_* |\theta(t) -\theta_{0, \ell}|_{L^2(\Omega)}^2 \bigr)
\\[1ex]
\ds +2 \bigl( |\eta_{0, \ell} -\eta_0|_{L^2(\Omega)}^2 +A_* |\theta_{0, \ell} -\theta_0|_{L^2(\Omega)}^2 \bigr), \mbox{ \ for any $ t \in (0, \infty) $ and $ \ell \in \N $,}
\end{array}
\end{equation*}
we infer from (\ref{kenV-30}) that
\begin{equation*}
\begin{array}{c}
\ds \limsup_{t \downarrow 0} \, \bigl( |\eta(t) -\eta_0|_{L^2(\Omega)}^2 +A_* |\theta(t) -\theta_0|_{L^2(\Omega)}^2 \bigr)
\\[1ex]
\ds \leq 4 \bigl( |\eta_{0, \ell} -\eta_0|_{L^2(\Omega)}^2 +A_* |\theta_{0, \ell} -\theta_0|_{L^2(\Omega)}^2 \bigr), \mbox{ \ for any $ \ell \in \N $.}
\end{array}
\end{equation*}
{
By (\ref{kenV-01}), (\ref{kenV-29}) and (\ref{kenV-06})-(\ref{kenV-07}), the above inequality implies (\hyperlink{2-b}{$\sharp$2-b}).
}
\medskip

Next, given any bounded open interval $ I \subset \hspace{-0.25ex} \subset (0, \infty) $, we take a sequence $ \{ \eta_n \}_{n = 1}^\infty \subset C^\infty(\overline{I \times \Omega}) $, such that $ \eta_n \to \eta $ in $ L^2(I; H^1(\Omega)) $ as $ n \to \infty $. We choose  $ \omega = \overline{\eta}_n -\eta_n $ in  (\ref{kenV-09}). Then, having in mind (\ref{kenV-34})-{(\ref{1steq-meas})}, and applying {Lemma \ref{Lem.A_last}} with $ \beta = 1 $, $ \{ \beta_n \}_{n = 1}^\infty = \{ 1 \} $, $ v = \theta $, $ \{ v_n \}_{n = 1}^\infty = \{ \underline{\theta}_n \}_{n = 1}^\infty $, $ \varrho = 0 $ and $ \{ \varrho_n \}_{n = 1}^\infty = \{ (\overline{\eta}_n -\eta_n) \alpha'(\overline{\eta}_n) \}_{n = 1}^\infty $, one can see that
\begin{equation}\label{kenV-36}
\begin{array}{rl}
\multicolumn{2}{l}{\ds\hspace{-2ex} \int_I |\nabla \eta(t)|_{L^2(\Omega)^N}^2 \, dt \leq \liminf_{n \to \infty} \int_I |\nabla \overline{\eta}_n(t)|_{L^2(\Omega)^N}^2 \, dt
\leq \limsup_{n \to \infty} \int_I |\nabla \overline{\eta}_n(t)|_{L^2(\Omega)^N}^2 \, dt
}
\\[2ex]
\leq & \ds \lim_{n \to \infty} \left[ \int_I |\nabla \eta_n(t)|_{L^2(\Omega)^N}^2 \, dt -2\int_I \bigl( (\widehat{\eta}_n)_t(t) +g(\overline{\eta}_n(t)), (\overline{\eta}_n -\eta_n)(t) \bigr)_{L^2(\Omega)} \, dt \right.
\\[2ex]
& \ds \quad  \left. -2\int_I \int_\Omega (\overline{\eta}_n -\eta_n)(t) \alpha'(\overline{\eta}_n(t)) {|\nabla \underline{\theta}_n(t)|_{\nu_n}} \, dx dt \right]
= \int_I |\nabla \eta(t)|_{L^2(\Omega)^N}^2\, dt.
\end{array}
\end{equation}
By (\ref{kenV-06})-(\ref{kenV-07}), (\ref{kenV-18}), (\ref{kenV-36}) and the uniform convexities of the $ L^2 $-type topologies, we have
{
\begin{equation}\label{kenV-36_1}
\left\{ \parbox{8.5cm}{
$ \overline{\eta}_n \to \eta $ in  $ L^2(I; H^1(\Omega)) $,
\\[1ex]
$ \ds \int_I \mathscr{F}_{\nu_n}(\overline{\eta}_n(t), \overline{\theta}_n(t)) \, dt \to \int_I \mathscr{F}(\eta(t), \theta(t)) \, dt $,
} \right. \mbox{as $ n \to \infty $}
\end{equation}
}
and,  by (\ref{kenV-02}),
\begin{equation}\label{kenV-37}
\begin{array}{rl}
\multicolumn{2}{l}{\ds \hspace{-4ex}\left| \int_I \mathscr{F}_{\nu_n}(\underline{\eta}_n(t), \underline{\theta}_n(t)) \, dt -\int_I \mathscr{F}(\eta(t), \theta(t)) \, dt \right|}
\\[2ex]
\leq & \ds \left| \int_I \mathscr{F}_{\nu_n}(\overline{\eta}_n(t), \overline{\theta}_n(t)) \, dt -\int_I \mathscr{F}_{\nu_n}(\underline{\eta}_n(t), \underline{\theta}_n(t)) \, dt \right|
\\[2ex]
& \ds \hspace{17.8ex} +\left| \int_I \mathscr{F}_{\nu_n}(\overline{\eta}_n(t), \overline{\theta}_n(t)) \, dt -\int_I \mathscr{F}(\eta(t), \theta(t)) \, dt \right|
\\[2ex]
\leq & \ds 2 h_n \mathscr{F}_{\nu_n}(\eta_{0, n}, \theta_{0, n}) +\left| \int_I \mathscr{F}_{\nu_n}(\overline{\eta}_n(t), \overline{\theta}_n(t)) \, dt -\int_I \mathscr{F}(\eta(t), \theta(t)) \, dt \right|
\\[2ex]
\to & 0, \mbox{ as $ n \to \infty $, \ for any bounded open interval $ I \subset \hspace{-0.25ex} \subset (0, \infty) $.}
\end{array}
\end{equation}

Given now any bounded open set $ A \subset (0, \infty) $, we denote by $ \mathcal{I}_A $ the (at most countable) class of pairwise-disjoint open intervals, such that \ $ \bigcup_{\tilde{I} \in \mathcal{I}_A} \tilde{I} $. Here, from (\ref{kenV-37}), it is deduced that
\begin{equation*}
\begin{array}{c}
\ds \sum_{\tilde{I} \in \tilde{\mathcal{I}}} \int_{\tilde{I}} \mathscr{F}(\eta(t), \theta(t)) \, dt \leq \liminf_{n \to \infty} \int_{A} \mathscr{F}_{\nu_n}(\underline{\eta}_n(t), \underline{\theta}_n(t)) \, dt
\\[2ex]
\mbox{for any finite subclass $ \tilde{\mathcal{I}} \subset \mathcal{I}_A $,}
\end{array}
\end{equation*}
and therefore,
\begin{equation}\label{kenV-38}
\begin{array}{c}
\ds \int_A \mathscr{F}(\eta(t), \theta(t)) \, dt \leq \liminf_{n \to \infty} \int_A \mathscr{F}_{\nu_n}(\underline{\eta}_n(t), \underline{\theta}_n(t)) \, dt
\\[2ex]
\mbox{for any bounded open set $ A \subset \hspace{-0.25ex} \subset (0, \infty) $.}
\end{array}
\end{equation}
As an application of \cite[Proposition 1.80]{AFP} for (\ref{kenV-36_1}), (\ref{kenV-37}) and (\ref{kenV-38}), we obtain
\begin{equation}\label{kenV-39}
\mathscr{F}_{\nu_n}(\underline{\eta}_n, \underline{\theta}_n) \to \mathscr{F}(\eta, \theta) \mbox{ in $ \mathcal{M}_{\rm loc}((0, \infty)) $, as $ n \to \infty $.}
\end{equation}

Therefore, by (\ref{kenV-08}),
\begin{equation}\label{kenV-40}
\mathscr{F}(\eta(t), \theta(t)) = \mathscr{J}_*(t), \mbox{ \ a.e. $ t \in (0, \infty) $.}
\end{equation}
Finally, having (\ref{kenV-08}) in mind and letting $ n \to \infty $ in (\ref{kenV-100}) yields
\begin{equation}\label{kenV-41}
\begin{array}{c}
\ds \frac{1}{2} \int_s^t |\eta_t(\tau)|_{L^2(\Omega)}^2 \, d \tau +\int_s^t |{\textstyle \sqrt{\alpha_0(\eta(\tau))}} \theta_t(\tau)|_{L^2(\Omega)}^2 \, d \tau +\mathscr{J}_*(t) \leq \mathscr{J}_*(s),
\\[2ex]
\mbox{for a.e. $ 0 < s < t < \infty $.}
\end{array}
\end{equation}
{
From \eqref{kenV-06} and (\hyperlink{2-a}{$\sharp$2-a}), we observe that
}
\begin{center}
$ \mathscr{F}(\eta, \theta) \in L_{\rm loc}^1([0, \infty)) $ $\cap L_{\rm loc}^\infty((0, \infty)) $,
\end{center}
which, together with (\ref{kenV-40})-(\ref{kenV-41}) leads to (\hyperlink{2-c}{$\sharp$2-c}). \hfill $ \Box $

\section{Proof of Main Theorem 2}
\ \ \vspace{-3ex}

{
\begin{rem}\label{Rem.extra01}
\begin{em}
In what follows, we assign the left-continuous expression of $ t \in (0, \infty) \mapsto \mathscr{F}(\eta(t), \theta(t)) $ to the function $ \mathscr{J}_* $ in (\hyperlink{S4}{S4}) of Definition \ref{sol}.
Then, due to the nonincreasing property of $ \mathscr{J}_* $,  the condition ``for a.e. $ 0 < s < t < \infty $'' in (\ref{kenV-41}) can be rephrased as ``for all $ 0 < s \leq t < \infty $''. Moreover, taking into account Remark \ref{Rem.auxCor1} and (\hyperlink{(S0)}{S0}) in Definition \ref{sol}, one can deduce from (\ref{kenV-41}) that
\begin{equation}\label{ene-inq1}
\begin{array}{c}
\ds \frac{1}{2} \int_s^t |\eta_t(\tau)|_{L^2(\Omega)}^2 \, d \tau +\int_s^t |{\textstyle \sqrt{\alpha_0(\eta(\tau))}} \theta_t(\tau)|_{L^2(\Omega)}^2 \, d \tau +\mathscr{F}(\eta(t), \theta(t)) \leq \mathscr{J}_*(s),
\\[2ex]
\mbox{for all $ 0 < s \leq t < \infty $.}
\end{array}
\end{equation}
\end{em}
\end{rem}
}
\bigskip

\noindent
\textbf{Proof of Main Theorem \ref{mainTh2}.}
Let $ [\eta, \theta] \in C([0, \infty); L^2(\Omega))^2 \cap W_{\rm loc}^{1, 2}((0, \infty); L^2(\Omega))^2 $ be an energy-dissipative solution to (\hyperlink{S}{S}). Then, from (\ref{ene-inq1}) and Remark \ref{Rem.auxCor1}, it is observed that
\begin{equation}\label{kenVI-01}
\left\{ \hspace{-2.5ex} \parbox{14cm}{
\vspace{-2ex}
\begin{itemize}
\item $ [\eta_t, \theta_t] \in L^2(1, \infty; L^2(\Omega))^2 $,
\item $ \{ [\eta(t), \theta(t)] \, | \, t \geq 1 \} \subset F_1 := \left\{ \begin{array}{l|l}
[\tilde{\eta}, \tilde{\theta}] \in D_*(\theta_0) & \mathscr{F}(\tilde{\eta}, \tilde{\theta}) \leq \mathscr{J}_*(1)
\end{array} \right\} $,
\item $ F_1 $ is compact in $ L^2(\Omega)^2 $.
\vspace{-2ex}
\end{itemize}
} \right.
\end{equation}
Therefore, we can find a pair $ [\eta_\infty, \theta_\infty] \in L^2(\Omega)^2 $ and a sequence of times $1 \leq t_1 < t_2 < t_3 < \cdots < t_n  \uparrow \infty$ as $ n \to \infty$
such that
\begin{equation}\label{kenVI-02}
[\eta(t_n), \theta(t_n)] \to [\eta_\infty, \theta_\infty] \mbox{ in $ L^2(\Omega)^2 $, as $ n \to \infty $.}
\end{equation}
This implies that $ \omega(\eta, \theta) \ne \emptyset $. Also, the compactness of $ \omega(\eta, \theta) $ is obtained from the compactness of $ F_1 $ and the fact that
\begin{equation*}
\omega(\eta, \theta) = \bigcap_{s \geq 0} \overline{\{ [\eta(t), \theta(t)] \, | \, t \geq s \}} \subset \overline{\{ [\eta(t), \theta(t)] \, | \, t \geq 1 \}} \subset F_1.
\end{equation*}
Thus, (\hyperlink{(O)}{O}) holds.

\medskip
Next, we verify  (\,\hyperlink{(I)}{I}\,). We take any $ [\eta_\infty, \theta_\infty] \in \omega(\eta, \theta) $. Then, there exists a sequence of times $1 \leq t_1 < t_2 < t_3 < \cdots < t_n  \uparrow \infty$ as $ n \to \infty$
such that
(\ref{kenVI-02}) holds. Therefore,  item (\hyperlink{i-a}{i-a}) is a straightforward consequence of {(\hyperlink{(S0)}{S0})}, (\ref{kenVI-01}) and (\ref{kenVI-02}).

\medskip In the meantime, it follows from (\ref{kenVI-01}) that
\begin{equation}\label{kenVI-03}
\left\{ \hspace{-2.5ex} \parbox{12cm}{
\vspace{-1ex}
\begin{itemize}
\item $ \{ \eta_n \}_{n = 1}^\infty := \{ \eta({}\cdot +\,t_n) \}_{n = 1}^\infty $ is bounded in $ W^{1, 2}(0, 1; L^2(\Omega)) \cap L^\infty(0, 1; H^1(\Omega)) $;
\vspace{-1ex}
\item $ \{ \theta_n \}_{n = 1}^\infty := \{ \theta({}\cdot + \, t_n) \}_{n = 1}^\infty $ is bounded in $ W^{1, 2}(0, 1; L^2(\Omega)) $, and $ \{ |D \theta_n({}\cdot{})|(\Omega) \}_{n = 1}^\infty $ is bounded in $ L^\infty(0, 1) $;
\vspace{-1ex}
\item $ \{ [\eta_n(t), \theta_n(t)] \, | \, t \in [0, 1], ~ n \in \N \} \subset D_*(\theta_0) $.
\vspace{-1ex}
\end{itemize}
} \right.
\end{equation}
Owing to (\ref{kenVI-01}) and (\ref{kenVI-03}) { and the compactness theories as in \cite{AFP,Simon},} we  infer that
\begin{equation}\label{kenVI-06}
\parbox{12cm}{
\begin{tabular}{rl}
$ (\eta_n)_t \to 0 $ & \hspace{-1ex} and \ $ (\theta_n)_t \to 0 $ \ in $ L^2(0, 1; L^2(\Omega)) $, \ as $ n \to \infty $,
\end{tabular}
}
\end{equation}
\begin{equation}\label{kenVI-04}
\parbox{12cm}{
\begin{tabular}{rl}
$ \eta_n \to \eta_\infty $ & in $ W^{1, 2}(0, 1; L^2(\Omega)) $, weakly-$*$ in $ L^\infty(0, 1; H^1(\Omega)) $,
\\[1ex]
& and weakly-$*$ in $ L^\infty((0, 1) \times \Omega) $, \ as $ n \to \infty $;
\end{tabular}
}
\end{equation}
\begin{equation}\label{kenVI-05}
\parbox{12cm}{
\begin{tabular}{rl}
$ \theta_n \to \theta_\infty $ & in $ W^{1, 2}(0, 1; L^2(\Omega)) $, \ as $ n \to \infty $;
\end{tabular}
}
\end{equation}
\begin{equation}\label{kenVI-07}
\parbox{12cm}{
\begin{tabular}{rl}
$ \theta_n(t) \to \theta_\infty $ & weakly-$ * $ in $ BV(\Omega) $, for any $ t \in (0, 1) $, \ as $ n \to \infty $,
\end{tabular}
}
\end{equation}
{
by taking subsequences (not relabeled) if necessary.}
Additionally, due to (\hyperlink{(S1)}{S1})-(\hyperlink{(S2)}{S2}) in Definition \ref{sol}, the sequence $ \{ [\eta_n, \theta_n] \}_{n = 1}^\infty $ satisfies
\begin{equation}\label{kenVI-08}
\begin{array}{c}
\ds \int_0^1 \bigl( (\eta_n)_t(t) +g(\eta_n(t)), w \bigr)_{L^2(\Omega)} \, dt +\int_0^1 (\nabla \eta_n(t), \nabla w)_{L^2(\Omega)^N} \, dt
\\[1.5ex]
\ds +\int_0^1 \int_\Omega  d \bigl[ w \alpha'(\eta(t)) |D \theta(t)| \bigr] \, dt = 0,
\\[2ex]
\mbox{ \ for any $ w \in H^1(\Omega) \cap L^\infty(\Omega) $ and any $ n \in \N $,}
\end{array}
\end{equation}
and
\begin{equation}\label{kenVI-09}
\begin{array}{c}
\ds \int_0^1 \bigl( \alpha_0(\eta_n(t)) (\theta_n)_t(t), \theta_n(t) \bigr)_{L^2(\Omega)} \, dt +\int_0^1 \Phi_0(\alpha(\eta_n(t)); \theta_n(t)) \, dt
\\[2ex]
\leq \Phi_0(\alpha(\eta_n(t)); 0)  =0 , \mbox{ \ for any $ n \in \N $.}
\end{array}
\end{equation}

Taking into account (\ref{kenVI-01}), {(\ref{kenVI-06})-(\ref{kenVI-07})} and (\ref{kenVI-09}) and Lemma \ref{auxLem1}, we obtain
\begin{eqnarray}
0 & \leq & \int_0^1 \int_\Omega d \bigl[ \alpha(\eta(t)) |D \theta(t)| \bigr] dt \leq \liminf_{n \to \infty} \int_0^1 \int_\Omega d \bigl[ \alpha(\eta_n(t)) |D \theta_n(t)| \bigr] \, dt
\nonumber
\\
& \leq & \limsup_{n \to \infty} \int_0^1 \Phi_0(\alpha(\eta_n(t)), \theta_n(t)) \, dt
\label{kenVI-10}
\\
& \leq & -\lim_{n \to \infty} \int_0^1 \bigl( \alpha_0(\eta_n(t)) (\theta_n)_t(t), \theta_n(t) \bigr)_{L^2(\Omega)} \, dt
= 0.
\nonumber
\end{eqnarray}
By (\hyperlink{(H4)}{H4}), the above inequality implies the item (\hyperlink{i-c}{i-c}).
\medskip

Finally, with (\ref{kenVI-06})-(\ref{kenVI-07}) and (\ref{kenVI-10}) in mind, given any $ w \in H^1(\Omega) \cap L^\infty(\Omega) $, we apply Lemma \ref{auxLem2} with $ I = (0, 1) $, $ \beta = \alpha(\eta_\infty) $, $ \{ \beta_n \}_{n = 1}^\infty = \{ \alpha(\eta_n) \}_{n = 1}^\infty $, $ v =\theta_\infty $, $ \{ v_n \}_{n = 1}^\infty =  \{ \theta_n \}_{n = 1}^\infty $, $ \varrho = w \alpha'(\eta_\infty) $ and $ \{ \varrho_n \}_{n = 1}^\infty = \{ w \alpha'(\eta_n) \}_{n = 1}^\infty $. Then, we infer that
\begin{equation}\label{kenVI-11}
\begin{array}{c}
\ds \int_0^1 \int_\Omega d \bigl[ w \alpha'(\eta_n(t)) |D \theta_n(t)| \bigr] \, dt \to \int_\Omega d \bigl[ w \alpha'(\eta_\infty) |D \theta_\infty| \bigr] = 0,
\\[2ex]
\mbox{as $ n \to \infty $, \ for any $ w \in H^1(\Omega) \cap L^\infty(\Omega) $.}
\end{array}
\end{equation}
With (\ref{kenVI-06})-(\ref{kenVI-07}) and (\ref{kenVI-11}) in mind, letting $ n \to \infty $ in (\ref{kenVI-08}) yields
\begin{equation*}
(\nabla \eta_\infty, \nabla w)_{L^2(\Omega)^N} +(g(\eta_\infty), w)_{L^2(\Omega)} = 0, \mbox{ \ for any $ w \in H^1(\Omega) \cap L^\infty(\Omega) $.}
\end{equation*}
{Thus, we conclude (\hyperlink{i-b}{i-b}).
\hfill $ \Box $}


\label{page:e}
\end{document}